\documentclass[final,onefignum,onetabnum]{siamart171218}

\usepackage{lipsum}
\usepackage{amsfonts}
\usepackage{graphicx}
\usepackage{epstopdf}
\usepackage{algorithmic}
\usepackage{multirow}
\usepackage{epstopdf}
\usepackage{mathtools}
\usepackage{booktabs}
\usepackage{amssymb,relsize}
\usepackage{amsmath}
\usepackage[normalem]{ulem}
\usepackage[caption=false]{subfig}
\ifpdf
\DeclareGraphicsExtensions{.eps,.pdf,.png,.jpg}
\else
\DeclareGraphicsExtensions{.eps}
\fi

\DeclarePairedDelimiter\ceil{\lceil}{\rceil}

\DeclarePairedDelimiter\abs{\lvert}{\rvert}%

\makeatletter
\let\oldabs\abs
\def\abs{\@ifstar{\oldabs}{\oldabs*}}

\newsiamremark{remark}{Remark}
\newsiamremark{hypothesis}{Hypothesis}
\crefname{hypothesis}{Hypothesis}{Hypotheses}
\newsiamthm{claim}{Claim}

\headers{A Novel Partitioning Method for Block Cimmino}
{F.~Sukru~Torun, Murat~Manguoglu, Cevdet~Aykanat}

\title{A Novel Partitioning Method for Accelerating the Block Cimmino Algorithm  \thanks{Submitted to the editors January 3, 2018.
		\funding{The work of the first author was supported by the Scientific and Technological Research Council of Turkey (TUBITAK), under the program BIDEB-2211. 
			The second author would like to thank Alexander von Humboldt
			Foundation for support of a research stay at TU-Berlin.
}}}

\author{F.~Sukru~Torun
	\footnotemark{}
	\and
	Murat~Manguoglu 
	\footnotemark{}
	\footnotemark{}    
	\and
	Cevdet~Aykanat
	\footnotemark{}     
}

\usepackage{amsopn}

\ifpdf
\hypersetup{
  pdftitle={A Novel Partitioning Method for Accelerating the Block Cimmino Algorithm},
  pdfauthor={F.~Sukru~Torun, Murat~Manguoglu, and Cevdet~Aykanat}
}
\fi

\begin{document}
\maketitle

\renewcommand{\thefootnote}{\fnsymbol{footnote}}
\footnotetext[2]{Department of Computer Engineering, Ankara Y{\i}ld{\i}r{\i}m Beyaz{\i}t University, 06010, Ankara, Turkey  (\email{\mbox{ftorun@ybu.edu.tr}}).}
\footnotetext[3]{Institut f\"ur Mathematik, Technische Universit\"at Berlin, 10623 Berlin, Germany.}
\footnotetext[4]{Department of Computer Engineering, Middle East Technical University, 06800 Ankara, Turkey    (\email{\mbox{manguoglu@ceng.metu.edu.tr}}).}
\footnotetext[5]{Department of Computer Engineering, Bilkent University, 06800 Ankara, Turkey    (\email{aykanat@cs.bilkent.edu.tr}).}

\begin{abstract}
We propose a novel block-row partitioning method in order to improve the convergence rate of the block Cimmino algorithm for solving general sparse linear systems of equations. 
The convergence rate of the block Cimmino algorithm depends on the orthogonality among the block rows obtained by the partitioning method.
The proposed method takes numerical orthogonality among block rows into account by proposing a row inner-product graph model of the coefficient matrix.
In the graph partitioning formulation defined on this graph model, the partitioning objective of minimizing the cutsize directly corresponds to minimizing the sum of inter-block inner products
between block rows thus leading to an improvement in the eigenvalue spectrum of the iteration matrix.
This in turn leads to a significant reduction in the number of iterations required for convergence.
Extensive experiments conducted on a large set of matrices confirm  the validity of the proposed method against a state-of-the-art method.

\end{abstract}

\begin{keywords}
	row projection methods, block Cimmino algorithm, Krylov subspace methods, row inner-product graph, graph partitioning
\end{keywords}

\begin{AMS}
  65F10, 65F50, 05C50, 05C70
\end{AMS}

\section{Introduction}
\label{intro}

Row projection methods are a class of iterative linear system solvers
\cite{bramley1989row,bramley1992row,kamath1989projection} that are used for solving linear system of equations of the form 
\begin{equation}
\label{eq:systemAx}
Ax=f,
\end{equation} 
where $A$ is an $n \times n$ sparse nonsymmetric nonsingular matrix, $x$ and $f$ are column vectors of size $n$.
In these methods, the solution is computed through successive projections onto rows of $A$.
There are mainly two major variations known as Kacmarz \cite{kaczmarz1937angenaherte} and Cimmino \cite{cimmino1938calcolo}.  
Kaczmarz obtains the solution through a product of orthogonal projections whereas Cimmino reaches the solution through a sum of orthogonal projections.
Cimmino is known to be more amenable to parallelism than Kaczmarz \cite{bramley1992row}.
However, Kaczmarz can be still parallelized via block Kaczmarz~\cite{kamath1989projection}, CARP~\cite{gordon2005} or multi-coloring~\cite{galgon2015parallel}.
The required number of iterations for Cimmino algorithm, however, could be quite large. 
One alternative variation is the block Cimmino~\cite{arioli1992block} which is a block row projection method.  
Iterative block Cimmino has been studied extensively in  \cite{arioli1992block,arioli1995block,bramley1992row,censor1988parallel,drummond2015partitioning,
elfving1980block,elfving2009properties}. 
A pseudo-direct version of block Cimmino based on the augmented block row projections is proposed in~\cite{duff2015augmented}.
However as in any other direct solver~\cite{MUMPS,bolukbasi2016multithreaded,li2005overview,schenk2004solving}, this approach also suffers from extensive memory requirement due to fill-in.

In the block Cimmino scheme, the linear system~\cref{eq:systemAx} is partitioned into $K$ block rows,
where $K \leq n $, as follows: 
\begin{equation} 
\label{eq:blockCimmino} 
\begin{pmatrix} 
A_1  \\
A_2  \\
\vdots  \\
A_K  \\
\end{pmatrix} 
x = 
\begin{pmatrix} 
f_1 \\
f_2 \\
\vdots \\
f_K \\
\end{pmatrix}.
\end{equation} 
\noindent In~\cref{eq:blockCimmino}, the coefficient matrix and right-hand side vector are partitioned conformably. Here $A_i$ is a submatrix of size $m_i \times n$ and $f_i$ is a subvector of size $m_i$, for $i=1,2,...,K$, where
\begin{equation}
n = \sum\limits_{i=1}^{K}m_i.
\end{equation}

The block Cimmino scheme is given in~\Cref{algo1}, where $A_i^+$ is the Moore-Penrose pseudoinverse of $A_i$ and it is defined as \begin{equation}
A_i^+ = A_i^T(A_iA_i^T)^{-1}.
\end{equation}
In~\Cref{algo1}, $A_i^+$ is used for the sake of the clarity of the notation and it is never computed explicitly.
In fact, at line 4 of~\Cref{algo1}, the minimum norm solution of an underdetermined linear least squares problem is obtained via the augmented system approach as discussed later.
In~\Cref{algo1}, $\delta_i$ vectors, which are of size $n$, can be computed independently in parallel without any communication, hence the block Cimmino algorithm is quite suitable for parallel computing platforms.  
At line 6, the solution is updated by the sum of  projections which is scaled by the relaxation parameter ($\omega$). 
In the parallel Cimmino algorithm, communication is required only for summing up the $\delta_i$s.

\begin{algorithm}
\caption{Block Cimmino method}
\label{algo1}
\begin{algorithmic}[1]
\STATE{Choose $x^{(0)}$}
\WHILE{$t=0,1,2,\ldots,$ until convergence}
\FOR{$i = 1,\ldots,K$}
\STATE{$\delta_i = A_i^+(f_i - A_ix^{(t)})$}
\ENDFOR
\STATE{$x^{(t+1)} = x^{(t)} + \omega \sum\limits_{i=1}^{K} \delta_i $ }
\ENDWHILE
\end{algorithmic}
\end{algorithm}

An iteration of  block Cimmino method can be reformulated as follows:
\begin{equation}
\label{eq:iterscheme}
\begin{tabular}{ll}   
$x^{(t+1)}  $ & $ =  x^{(t)} + \omega \sum\limits_{i=1}^{K} A_i^+ \left(  f_i - A_ix^{(t)} \right) $ \\
& $ = \left( I - \omega \sum\limits_{i=1}^{K} A_i^+A_i   \right)x^{(t)} + \omega \sum\limits_{i=1}^{K} A_i^+f_i $ \\
& $ = (I - \omega H) x^{(t)} +  \omega \sum\limits_{i=1}^{K} A_i^+f_i   $ \\
& $ = Qx^{(t)} + \omega \sum\limits_{i=1}^{K} A_i^+f_i   $,
\end{tabular} 
\end{equation}
where $Q$ is the iteration matrix for block Cimmino algorithm.
$\omega H = I- Q$ is the sum of $\mathcal{P_R}(A_i^T)$s (projections onto $A_i^T$) and it is defined by
\begin{equation}
\label{eq:H}
\begin{split}
\omega H & =  \omega \sum\limits_{i=1}^{K} \mathcal{P_R}(A_i^T) =   \omega \sum\limits_{i=1}^{K} A_i^+A_i.
\end{split}
\end{equation}

The projections in block Cimmino iterations can be calculated using several approaches, such as normal equations~\cite{golub2012matrix},  seminormal equations~\cite{golub1969,golub2012matrix}, QR factorization~\cite{golub2012matrix} and augmented system \cite{arioli1989augmented,golub2012matrix}.
The normal and seminormal equation approaches are not considered since they have the potential of introducing numerical difficulties that can be disastrous \cite{demmel1993improved,golub1969} in some cases when the problem is ill-conditioned.
Although the QR factorization is numerically more stable, it is computationally expensive.
Therefore we use the augmented system approach, which requires the solution of a sparse linear systems that can be done effectively by using a sparse direct solver~\cite{arioli1992block}.
Note that if submatrix $A_i$ is in a column overlapped block diagonal form, one could also use the algorithm in~\cite{torun2017parallel}.
However, this approach is not considered since we assume no structure for the coefficient matrix.

In the augmented system approach, we obtain the solution of~\cref{augmented} by solving the linear system~\cref{augmentedMatrix},
\begin{equation}
A_i \delta_i = r_i \hspace{1cm} (r_i=f_i - A_i x^{(t)})
\label{augmented}
\end{equation}
\begin{equation}
\begin{pmatrix}
    I & A_i^T  \\ 
    A_i &  0  \\ 
\end{pmatrix} 
\begin{pmatrix}
    \delta_i   \\ 
    \varsigma_i   \\ 
\end{pmatrix}
=
\begin{pmatrix}
    0   \\ 
    r_i   \\ 
\end{pmatrix}.
\label{augmentedMatrix}
\end{equation}
Hence, the solution of the augmented system gives $\delta_i$.

\subsection{The Conjugate Gradient acceleration of block Cimmino method}
\label{subsectioncgaccel}
Convergence rate of the block Cimmino algorithm is known to be slow~\cite{bramley1992row}.
In~\cite{bramley1992row}, the  Conjugate Gradient (CG) method is proposed to accelerate the row projection method. 
It is also reported that the CG accelerated Cimmino method competes favorably compared to classical preconditioned Generalized Minimum Residual (GMRES) and  Conjugate Gradient on Normal Equations (CGNE) for the solution of nonsymmetric linear systems that arise in an elliptic partial differential equation. 
On the other hand, it should be noted that the main motivation of block Cimmino algorithm is its amenability of parallelism \cite{drummond2015partitioning}.

The iteration scheme of the block Cimmino \cref{eq:iterscheme} gives
\begin{equation}
x^{(t+1)}  = (I-\omega H)x^{(t)} + \omega \sum\limits_{i=1}^{K} A_i^+f_i ,
\end{equation}
where the $H$ matrix is symmetric and positive definite according to~\cref{eq:H} if $A$ is square and full rank.
Hence, one can solve the following system using CG,
\begin{equation}
\label{eq:cg}
 \omega Hx = \omega\xi,
\end{equation}
where $\xi = \sum_{i=1}^{K} A_i^+f_i$ and $x$ is the solution vector of the system~\cref{eq:systemAx}. 
Note that since $\omega$ appears on both sides of~\cref{eq:cg}, it does not affect the convergence of CG.
\Cref{algo2} is the pseudocode for the CG accelerated block Cimmino method~\cite{zenadi2013methodes}, which is in fact the classical CG applied on the system given in~\cref{eq:cg}.
In the second line of the algorithm, the initial residual is computed in the same way as the first iteration of \cref{algo1}.
The matrix vector multiplications in the CG is expressed as the solution of $K$ independent underdetermined systems at line $5$ which can be done in parallel and  need to be summed to obtain $\psi^{(t)}$ by using an all-reduce operation.

\begin{algorithm}
\caption{Conjugate Gradient acceleration of block Cimmino method}
\label{algo2}
\begin{algorithmic}[1]
\STATE{Choose $x^{(0)}$}
\STATE{$r^{(0)}= \xi - \sum_{i=1}^{K} A_i^+A_i \; x^{(0)}$}
\STATE{$p^{(0)} = r^{(0)}$}
\WHILE{$t=0,1,2,\ldots,$ until convergence} \vspace{0.1cm}
\STATE{$\psi^{(t)} =  \sum_{i=1}^{K} A_i^+A_i \; p^{(t)}$}  \vspace{0.1cm}
\STATE{$ \alpha^{(t)} =   ({{{r}}^{(t)}}^T {{r}}^{(t)}) / ({p^{(t)}}^T \psi^{(t)}) $} \vspace{0.05cm}
\STATE{$ x^{(t+1)} = x^{(t)} +  \alpha^{(t)} p^{(t)}  $} \vspace{0.05cm}
\STATE{$ r^{(t+1)} = r^{(t)} -  \alpha^{(t)} \psi^{(t)}$} \vspace{0.05cm}
\STATE{$ \beta^{(t)} = ({r^{(t+1)}}^T r^{(t+1)})  /  ({{{r}}^{(t)}}^T {{r}}^{(t)}) $} \vspace{0.05cm}
\STATE{$  p^{(t+1)} = r^{(t+1)} + \beta^{(t)} p^{(t)} $}
\ENDWHILE
\end{algorithmic}
\end{algorithm}

\subsection{The effect of partitioning}   

The convergence rate of  CG accelerated block Cimmino algorithm is essentially the convergence rate of  CG applied on~\cref{eq:cg}.  A well-known upper bound on the convergence rate can be given in terms of the extreme eigenvalues ($\lambda_{min}$ and $\lambda_{max}$) of the coefficient matrix.  Let 
\begin{equation} 
\kappa = \frac{\lambda_{max}}{\lambda_{min}} 
\end{equation} 
be the 2-norm  condition number of $H$. Then,  as in \cite{golub2012matrix}, an upper bound on the convergence rate of the CG accelerated block Cimmino  can be  given by 
\begin{equation} 
\frac{||x^{(t)} -x^{*}||_{H}}{||x^{(0)} -x^{*}||_{H}} \leq 2   \left( \frac{\sqrt{\kappa} -1 } {\sqrt{\kappa} + 1} \right)^{t}  
\end{equation} 
where $x^{*}$  is the exact solution and  $||y||_{H}= y^{T}Hy$. Furthermore, it was shown that the convergence rate of CG  not only depends on the extreme eigenvalues but also on the separation between those extreme eigenvalues and interior eigenvalues~\cite{vanderSluis1986} as well the clustering of the internal eigenvalues~\cite{jennings1977influence}.    
In summary, the convergence rate of the CG accelerated block Cimmino depends on the extreme eigenvalues and the number of clusters  as well as the quality of the clustering.

Therefore, the partitioning of the coefficient matrix $A$ into block rows is crucial for improving the convergence rate of the CG accelerated block Cimmino algorithm, since it can improve the eigenvalue distribution of $H$. 
Note that the eigenvalues of $H$ is only affected by the block-row partitioning of the coefficient matrix $A$ and independent of any column ordering~\cite{drummond2015partitioning}.

Let the QR factorization of $A_i^T$ be defined as 
\begin{equation}
Q_iR_i = A_i^T,
\end{equation}
where the matrices $Q_i$ and $R_i$ have dimensions of $n \times m_i$ and $m_i \times m_i$, respectively.
Then, the $H$ matrix can be written as follows \cite{arioli1992block};
\begin{equation}
\label{eq:Iterat}
{\renewcommand{\arraystretch}{1.7}
\begin{tabular}{ll}   
$ H $ & $ =  \sum\limits_{i=1}^{K} A_i^T(A_iA_i^T)^{-1}A_i$ \\
$ $ & $ =  \sum\limits_{i=1}^{K} Q_iQ_i^T$ \\
$ $ & $ =  (Q_1, \ldots ,Q_K)(Q_1, \ldots ,Q_K)^T $. \\
\end{tabular} 
}
\end{equation}
Since the eigenvalue spectrum of matrix $(Q_1, \ldots ,Q_K)(Q_1, \ldots ,Q_K)^T$ is the same as the eigenvalue spectrum of matrix $(Q_1, \ldots ,Q_K)^T(Q_1, \ldots ,Q_K)$ \cite{golub1965calculating}, $H$ is similar to 
\begin{equation}
\begin{pmatrix}
{\renewcommand{\arraystretch}{2.0}
\begin{tabular}{llll}  
   $I_{m_1\times m_1}$ & ${Q_1}^T Q_2$ & $ \ldots $ & ${Q_1}^TQ_K$ \\
   ${Q_2}^T Q_1 $  & $I_{m_2\times m_2} $ & $ \ldots $ & ${Q_2}^TQ_K$ \\  
    $\sbox0{\dots}\makebox[\wd0]{\vdots} $  & $ \hdots $ & $ \ddots $ & $ \sbox0{\dots}\makebox[\wd0]{\vdots}$ \\  
   ${Q_K}^T Q_1 $  &  ${Q_K}^T Q_2 $ & $ \ldots $ & $I_{m_K \times m_K}$ \\  
\end{tabular}
}
\end{pmatrix},
\label{HQR}
\end{equation}
where the singular values of matrix ${Q_i}^TQ_j$ represent the principal angles between the subspaces spanned by the rows of $A_i$ and $A_j$ \cite{bjorck1973numerical}. 
Hence, the smaller the off-diagonals of the matrix~\cref{HQR}, the more eigenvalues of  $H$ will be clustered around one by the Gershgorin theorem \cite{gershgorin1931uber}.  
Therefore, the convergence rate of the block Cimmino method  highly depends on the orthogonality among $A_i$ blocks. 
If $A_i$ blocks are more orthogonal to each other, row inner products between blocks would be small and hence the eigenvalues will be clustered around one.  
Consequently,  CG is expected to converge in a fewer number of iterations.

In the literature, Cuthill-Mckee (CM)~\cite{cuthill1969reducing} based partitioning strategies~\cite{drummond2015partitioning,ruiz1992solution,zenadi2013methodes} are utilized to define block rows using CM level set information on the normal equations of the original matrix.
These strategies benefit from the level sets of CM for creating the desired number of block rows.
In CM, nodes in a level set have the same distance from the starting node and these nodes have neighbors only in the previous and the next level sets. 
Therefore permuted matrix based on the ordering of the level sets constitutes a block tridiagonal structure.
These strategies may suffer from not reaching the desired the number of block rows due to smaller number of level sets.
They also suffer from obtaining quite unbalanced partitions due to a relatively larger sizes of some level sets.
Numerical values are only considered when dropping based filtering strategy is used.
Although filtering small elements on normal equations before applying CM allows more freedom in partitioning by increasing the number of level sets, however, it does not hold the properties of the strict two-block partitioning~\cite{arioli1995block,ruiz1992solution,zenadi2013methodes}.
In addition, it is difficult to determine the best filtering threshold value a priori and find a common threshold which is ``good'' for all matrices. 

In recent studies~\cite{drummond2015partitioning,zenadi2013methodes},  a hypergraph partitioning method is used to find a ``good" block-row partitioning for the CG accelerated block Cimmino algorithm.
It is reported that it performs better than CM based methods (with or without filtering small elements).
In hypergraph partitioning, the partitioning objective of minimizing the cutsize corresponds to minimizing the number of linking columns among row blocks,  where a linking column refers to a column that contains nonzeros in more than one block row.
This in turn loosely relates to increasing the structural orthogonality~\cite{li2006miqr} among row blocks.
Here, two rows are considered to be structurally more orthogonal if they have fewer nonzeros in the same columns. This measure  depends on only nonzero counts and ignores the numerical values of nonzeros.

In this work, we propose a novel graph theoretical block-row partitioning method for the CG accelerated block Cimmino algorithm.
For this purpose, we introduce a row inner-product graph model of a given matrix $A$ and then the problem of finding a ``good" block-row partitioning is formulated as a graph partitioning problem on this graph model. The proposed method takes the numerical orthogonality between block rows of $A$ into account.  In the proposed method, the partitioning objective of minimizing the cutsize directly corresponds to minimizing the sum of inter-block inner products between block rows thus leading to an improvement in the eigenvalue spectrum of  $H$.

The validity of the proposed method is confirmed against two baseline methods by conducting experiments on a large set of matrices that arise in a variety of real life applications.
One of the two baseline methods is the state-of-the-art hypergraph partitioning method introduced in~\cite{drummond2015partitioning}.
We conduct experiments to study the effect of the partitioning on the eigenvalue spectrum of $H$.
We also conduct experiments to compare the performance of three methods in terms of the number of CG iterations and parallel CG time to solution.
The results of these experiments show that the proposed partitioning method is significantly better than both baseline methods in terms of all of these performance metrics.
Finally, we compare the preprocessing overheads of the methods and show that the proposed method incurs much less overhead than the hypergraph partitioning method, thus allowing a better amortization.

The rest of the paper is organized as follows. \Cref{section2} presents the proposed partitioning method and its implementation.
In \cref{section3}, we present and discuss the experimental results.
Finally, the paper is concluded with conclusions and directions for future research in \cref{section4}.


\section{The proposed partitioning method}
\label{section2}

In this section, we first describe the row inner-product graph model and then show that finding a ``good" block-row partitioning can be formulated as a graph partitioning problem on this graph model.
Finally, we give the implementation details for the construction and partitioning of the graph model.
We refer the reader to \cref{appendix}  for a short background on graph partitioning.

\subsection{Row inner-product graph model}

In the row-inner product graph \linebreak $\mathcal{G}_{\rm{RIP}}(A) = ( \mathcal{V}, \mathcal{E})$ of matrix $A$, vertices represent the rows of matrix $A$ and edges represent nonzero inner products between rows.
That is,  $\mathcal{V}$ contains  vertex $v_i $ for each row ${\rm{r}}_i $ of matrix $A$.
$\mathcal{E}$ contains an edge $(v_i,v_j)$ only if the inner product of rows ${\rm{r}}_i$ and ${\rm{r}}_j$ is nonzero.
That is,
\begin{equation}
\begin{aligned}
\mathcal{E}=\{(v_i,v_j) \; | \; {\rm{r}}_i{\rm{r}}_j^T   \; \neq \;  0 \}.\\
\end{aligned}
\end{equation}

\noindent 
Each vertex $v_i$ can be associated with a unit weight or a weight  that is equal to the number of nonzeros in row ${\rm{r}}_i$, that is
\begin{equation}
\label{eq16}
w(v_i)=1 \; \text{ or } \; w(v_i)=nnz({\rm{r}}_i),
\end{equation} 
respectively.
Each edge $(v_i,v_j)$  is associated with a cost  equal to the absolute value of the respective inner product.
That is,
\begin{equation}
\label{eq14}
cost(v_i,v_j) \; = \;  \abs{{\rm{r}}_i{\rm{r}}_j^T}  \; \text{ for all } \; (v_i,v_j) \in \mathcal{E}.
\end{equation}

If we prescale the rows of coefficient matrix $A$ such that each row has a unit \mbox{2-norm}, then the cost of edge  $(v_i,v_j)$ will correspond to the cosine of the angle between the pair of rows ${\rm{r}}_i$ and ${\rm{r}}_j$.
Therefore we prescale the matrix and the right-hand side vector in order to improve the effectiveness of the proposed graph model.
We note that the convergence of block Cimmino algorithm is independent of row scaling~\cite{elfving1980block,ruiz1992solution,zenadi2013methodes}.

This graph is topologically equivalent to the standard graph representation of the symmetric matrix resulting from the sparse matrix-matrix multiplication operation $C=AA^T$.
That is, the sparsity pattern of $C$ corresponds to the adjacency matrix representation of $\mathcal{G}_{\rm{RIP}}$.
Each nonzero $c_{ij}$ of the resulting matrix $C$ incurs an edge $(v_i,v_j)$.
Since each nonzero entry $c_{ij}$ of $C$ is computed as the inner product of row ${\rm{r}}_i$ and row ${\rm{r}}_j$, the absolute value of the nonzero $c_{ij}$ determines the cost of the respective edge $(v_i,v_j)$.
That is, since $c_{ij} = {\rm{r}}_i{\rm{r}}_j^T $, we have $cost(v_i,v_j) = \abs{c_{ij}}$.

\begin{figure}[hbtp]
  \begin{center}    
      \subfloat [Sample matrix $A$] {
       \includegraphics[scale=.45]{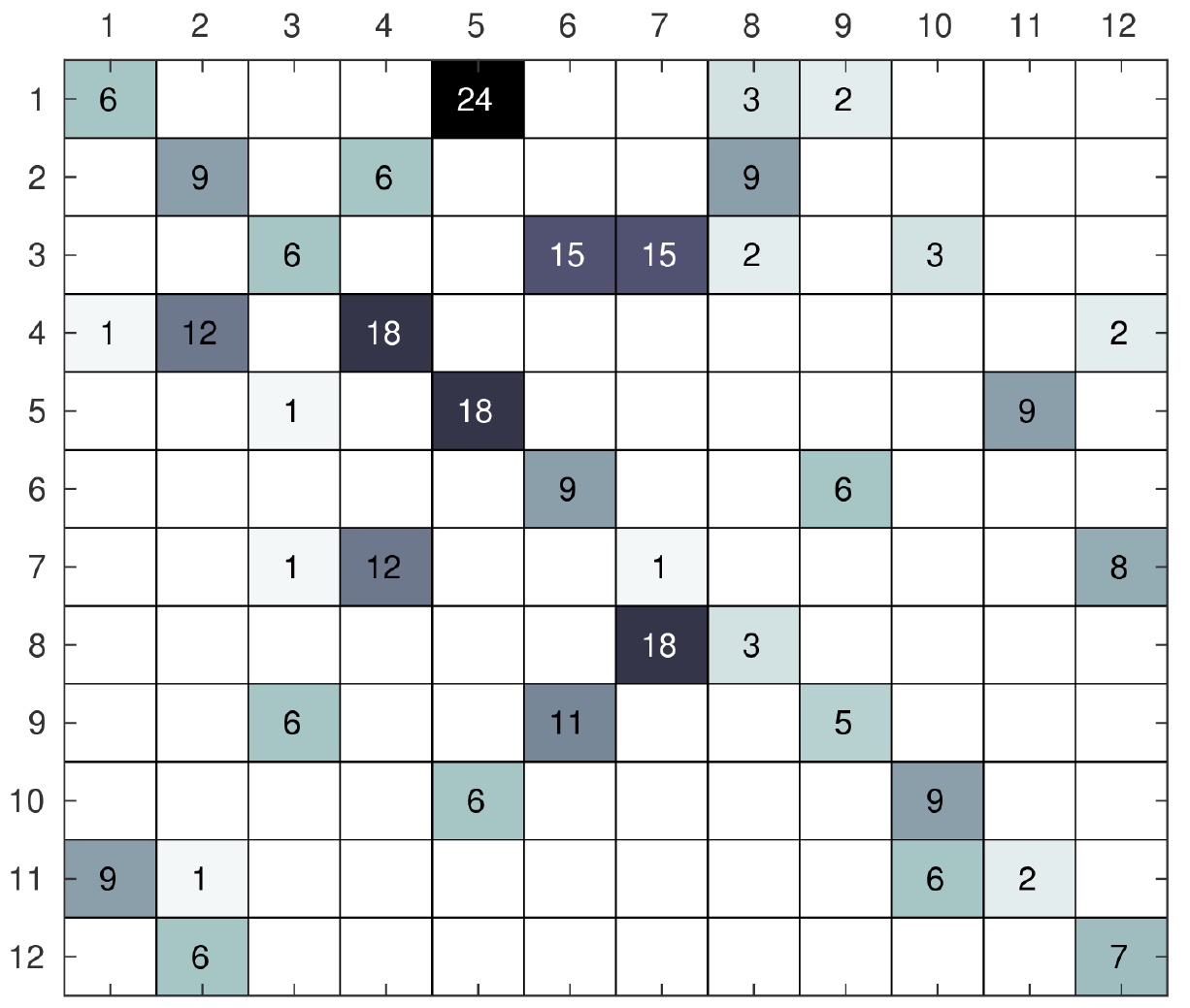}
        \label{graph11}
    }	\hspace{1em} 
      \subfloat [Row inner-product graph $\mathcal{G}_{\rm{RIP}}$ of A]{      
        \includegraphics[scale=.23]{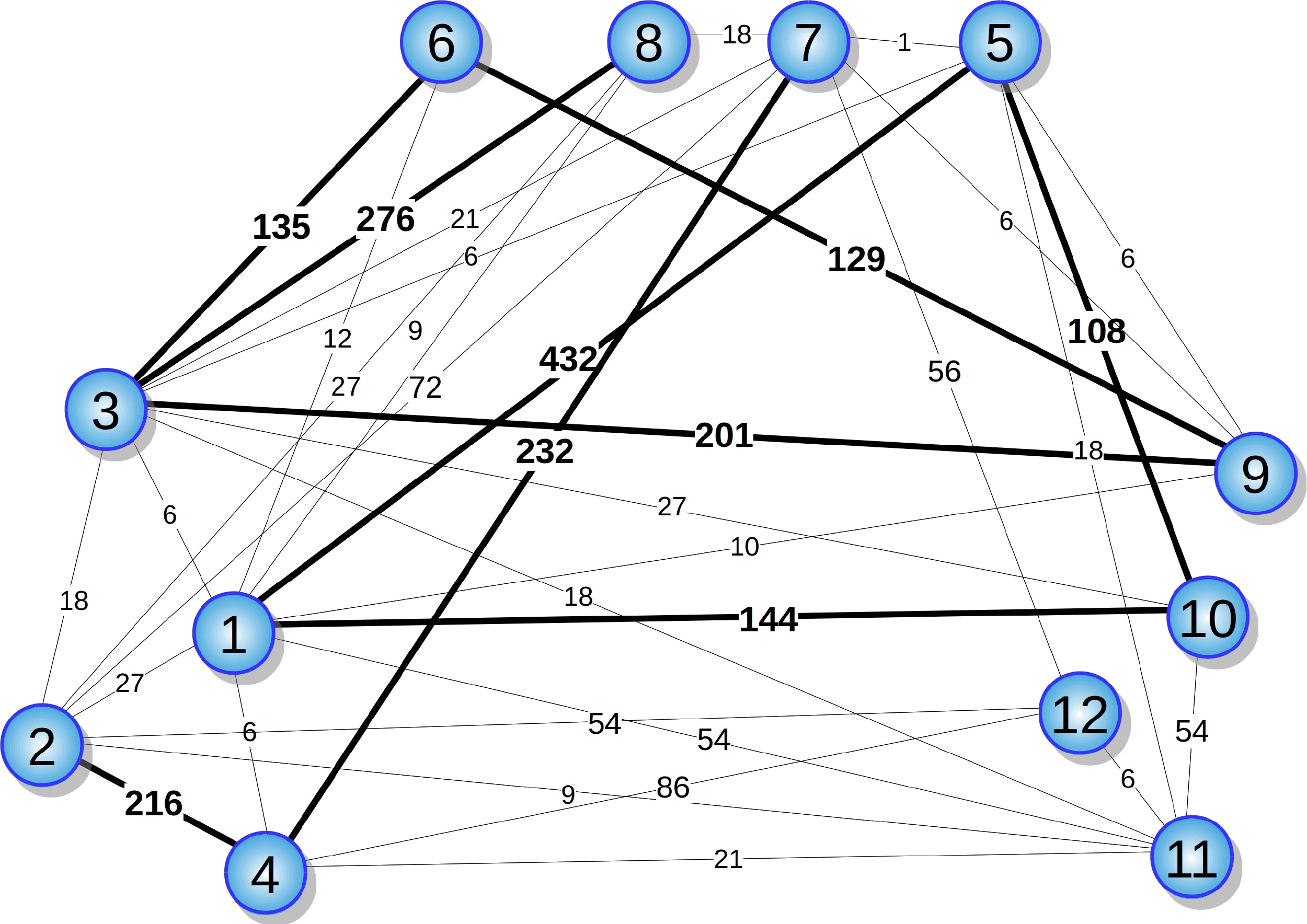}
        \label{graph12}
    }  \\    
 \end{center}  
	\caption{Row inner-product graph model.}
\label{fig:sample}
\vspace{-0.3cm}
\end{figure}

\Cref{graph11} shows a $12 \times 12$  sample sparse matrix that contains $38$ nonzeros.
Note that for the sake of clarity of presentation, rows of the sample matrix $A$ are not prescaled.
\Cref{graph12} depicts the proposed row inner-product graph $\mathcal{G}_{\rm{RIP}}$ for this sample matrix.
As seen in~\cref{graph12}, $\mathcal{G}_{\rm{RIP}}$ contains $12$ vertices each of which corresponds to a row and $35$ edges each of which corresponds to a nonzero row inner product.
For example, the inner product of rows ${\rm{r}}_2$ and ${\rm{r}}_4$ is nonzero where $  {\rm{r}}_2{\rm{r}}_4^T  = (9\! \times \!12) + (6\! \times \!18) = 216$ so that $\mathcal{G}_{\rm{RIP}}$ contains the edge $(v_2,v_4)$ with $cost(v_2,v_4)=216$.
In~\cref{graph12}, edges with cost larger than $100$ are shown with thick lines in order to make such high inner-product values more visible.

\Cref{aat} shows the resulting matrix $C$ of the sparse matrix-matrix multiplication $C= AA^T$. 
Only the off-diagonal nonzero entries together with their values are shown since the values of the diagonal entries do not affect $\mathcal{G}_{\rm{RIP}}$.
The cells that contain nonzeros larger than $100$ are shown with black background in order to make such high values more visible.
Comparison of~\cref{graph12,aat} shows that the topology of the standard graph model $\mathcal{G}(C)$ of matrix $C=AA^T$ is equivalent to the topology of $\mathcal{G}_{\rm{RIP}}$.
As also seen in~\cref{aat}, the values of the nonzero entries of matrix $C$ are equal to the costs of respective edges of $\mathcal{G}_{\rm{RIP}}$.
For example nonzero $c_{24}$ with a value $216$ incurs an edge $(v_2,v_4)$ with  $cost(v_2,v_4)=216$.

\subsection{Block-row partitioning via partitioning $\mathcal{G}_{\rm{RIP}}$}
A $K$-way partition $\Pi = \{ \mathcal{V}_1,\mathcal{V}_2,\ldots,\mathcal{V}_K \}$ of $\mathcal{G}_{\rm{RIP}}$ can be decoded as a partial permutation on the rows of $A$ to induce a permuted matrix $A^\Pi$, where 
\begin{equation}
A^\Pi = PA =  
\begin{bmatrix}
A_1^\Pi \\
\vdots \\
A^\Pi_k \\
\vdots \\
A^\Pi_K
\end{bmatrix} = \begin{bmatrix}
\mathcal{R}_1 \\
\vdots \\
\mathcal{R}_k \\
\vdots \\
\mathcal{R}_K
\end{bmatrix}.
\end{equation}

\noindent
Here, $P$ denotes the row permutation matrix which is defined by the $K$-way partition $\Pi$ as follows:
the rows associated with the vertices in $\mathcal{V}_{k+1}$ are ordered after the rows associated with the vertices in $\mathcal{V}_{k}$ for $k = 1,2,\ldots,K-1$.
That is, the block-row $\mathcal{R}_k$ contains the set of rows corresponding to the set of vertices in part $\mathcal{V}_k$ of partition $\Pi$, where  ordering of the rows within block row $\mathcal{R}_k$ is arbitrary for each $k=1,2,\ldots,K$. 
Note that we use the notation $\mathcal{R}_k$ to denote both $k^{th}$ block row $A^\Pi_k$ and the set of rows in $A^\Pi_k$.
Since the column permutation does not affect the convergence of block Cimmino algorithm~\cite{drummond2015partitioning}, the original column ordering of $A$ is maintained.

Consider a partition $\Pi$ of $\mathcal{G}_{\rm{RIP}}$.
The weight $W_k$ of a part $\mathcal{V}_k$ is either equal to the number of rows or number of nonzeros in block row $\mathcal{R}_k$ depending on the vertex weighting scheme used according to~\cref{eq16}.
That is,
\begin{equation}
W_k =  |\mathcal{R}_k| \; \text{  or  } \; W_k = \sum\limits_{r_i \in \mathcal{R}_k} nnz({\rm{r}}_i).
\end{equation}
\noindent
Therefore in partitioning $\mathcal{G}_{\rm{RIP}}$, the partitioning constraint of maintaining balance among part weights according to~\cref{partbal} corresponds to  maintaining balance on either the number of rows or the number of nonzeros among the block rows.

Consider a partition $\Pi$ of $\mathcal{G}_{\rm{RIP}}$.
A cut edge $(v_i,v_j)$ between parts  $\mathcal{V}_k$ and $\mathcal{V}_m$ represents a nonzero inter-block inner product ${\rm{r}}_i{\rm{r}}_j^T$  between block rows $\mathcal{R}_k$ and $\mathcal{R}_m$.
Therefore the cutsize of $\Pi$ (given in~\cref{partcut}) is equal to
\begin{equation}
\begin{aligned}
{\mathrm{cutsize}}(\Pi)  \triangleq \sum\limits_{ (v_i,v_j) \in \mathcal{E}_{cut}} cost(v_i,v_j) & =  \sum\limits_{1 \leq k < m \leq K } \sum\limits_{\substack{v_i \in \mathcal{V}_k \\ v_j \in \mathcal{V}_m}} cost(v_i,v_j) \\
& =   \sum\limits_{1 \leq k < m \leq K } \sum\limits_{\substack{{\rm{r}}_i \in \mathcal{R}_k \\ {\rm{r}}_j \in \mathcal{R}_m}}   \abs{{\rm{r}}_i{\rm{r}}_j^T},
 \end{aligned}
\end{equation}
\noindent
which corresponds to the total sum of the inter-block inner products ($\rm{ interIP}(\Pi)$).
So, in partitioning $\mathcal{G}_{\rm{RIP}}$,  the partitioning objective of minimizing the cutsize  corresponds to minimizing the sum of inter-block inner products between block rows.
Therefore this partitioning objective corresponds to making the block rows  numerically more orthogonal to each other.
This way, we expect this method to yield a faster convergence in the CG accelerated block Cimmino algorithm.

\begin{figure}[hbtp]
\vspace{-0.3cm}
  \begin{center}  
       \subfloat [] {
       \includegraphics[scale=.39]{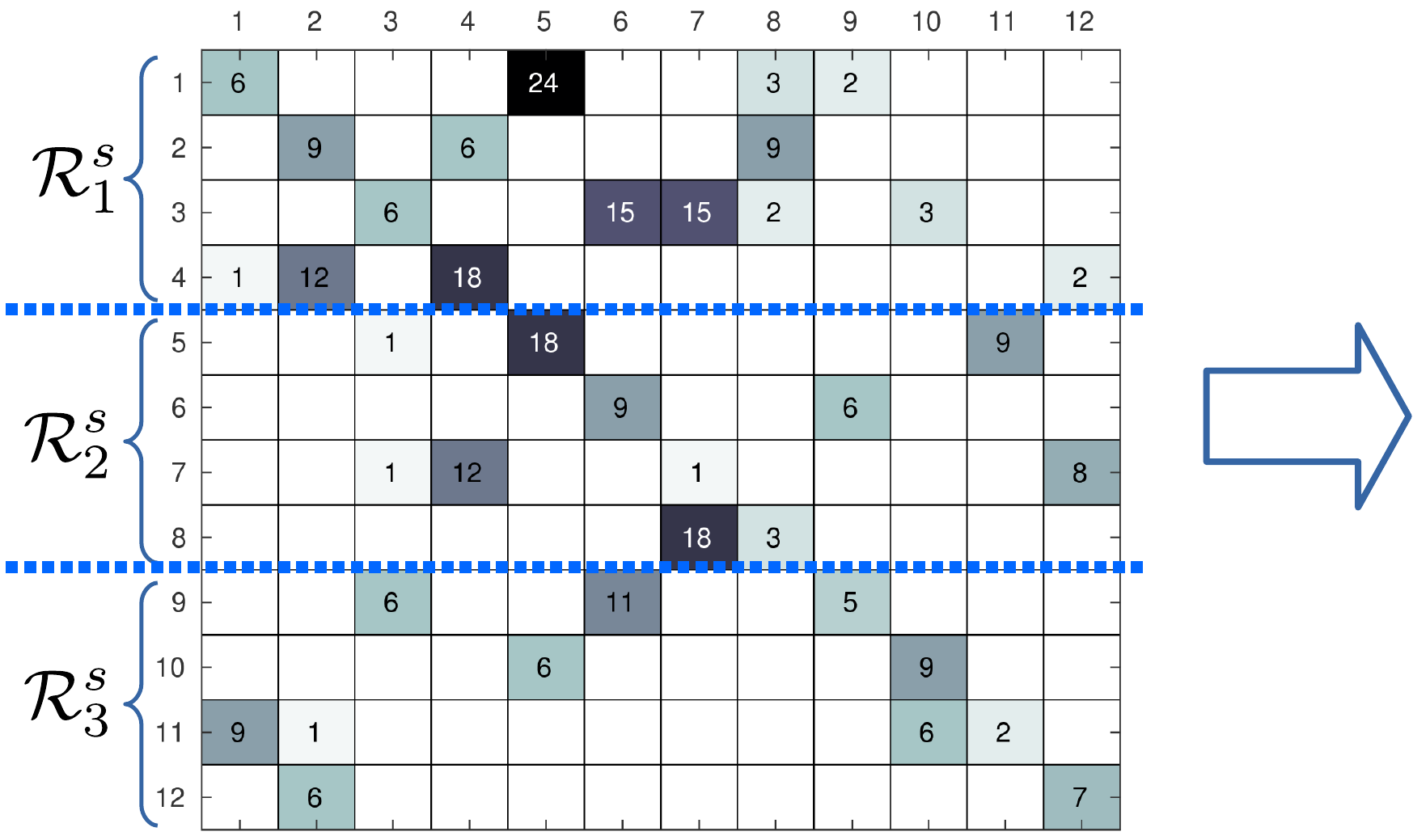}
        \label{unipart}
    }    
      \subfloat [] {
      \includegraphics[scale=.2]{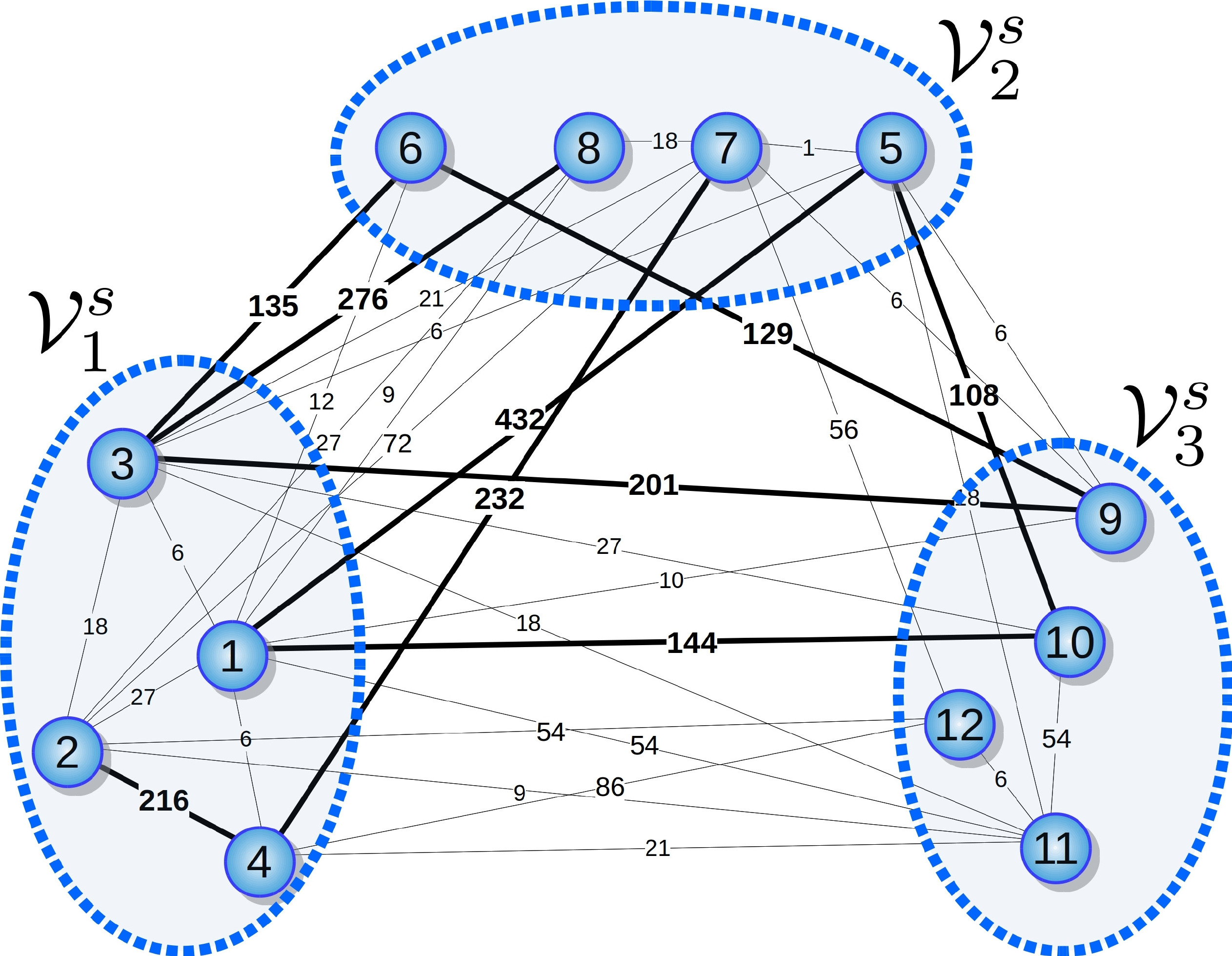}
        \label{uni}
    }	       	
 \end{center}  
	\caption{ (a) Straightforward 3-way row partition of $A$ and (b) 3-way partition $\Pi^s(\mathcal{V})$ of $\mathcal{G}_{\mathrm{{RIP}}}(A)$ induced by~\cref{unipart}.}
\label{graphpart1}
\end{figure}

\begin{figure}[hbtp]
  \begin{center}  
      \subfloat []{
        \includegraphics[scale=.2]{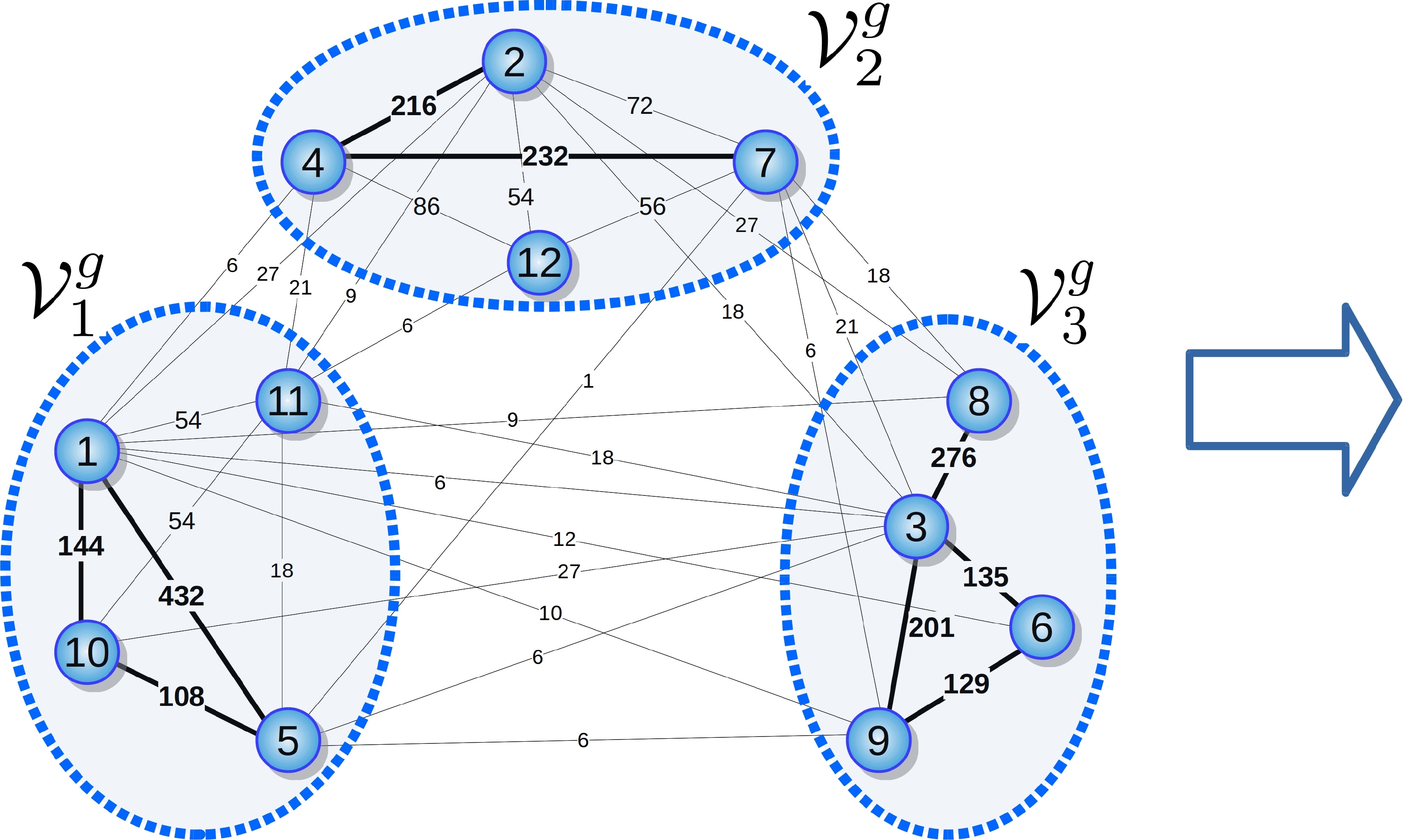}
        \label{met}
    }
     \subfloat []{
        \includegraphics[scale=.39]{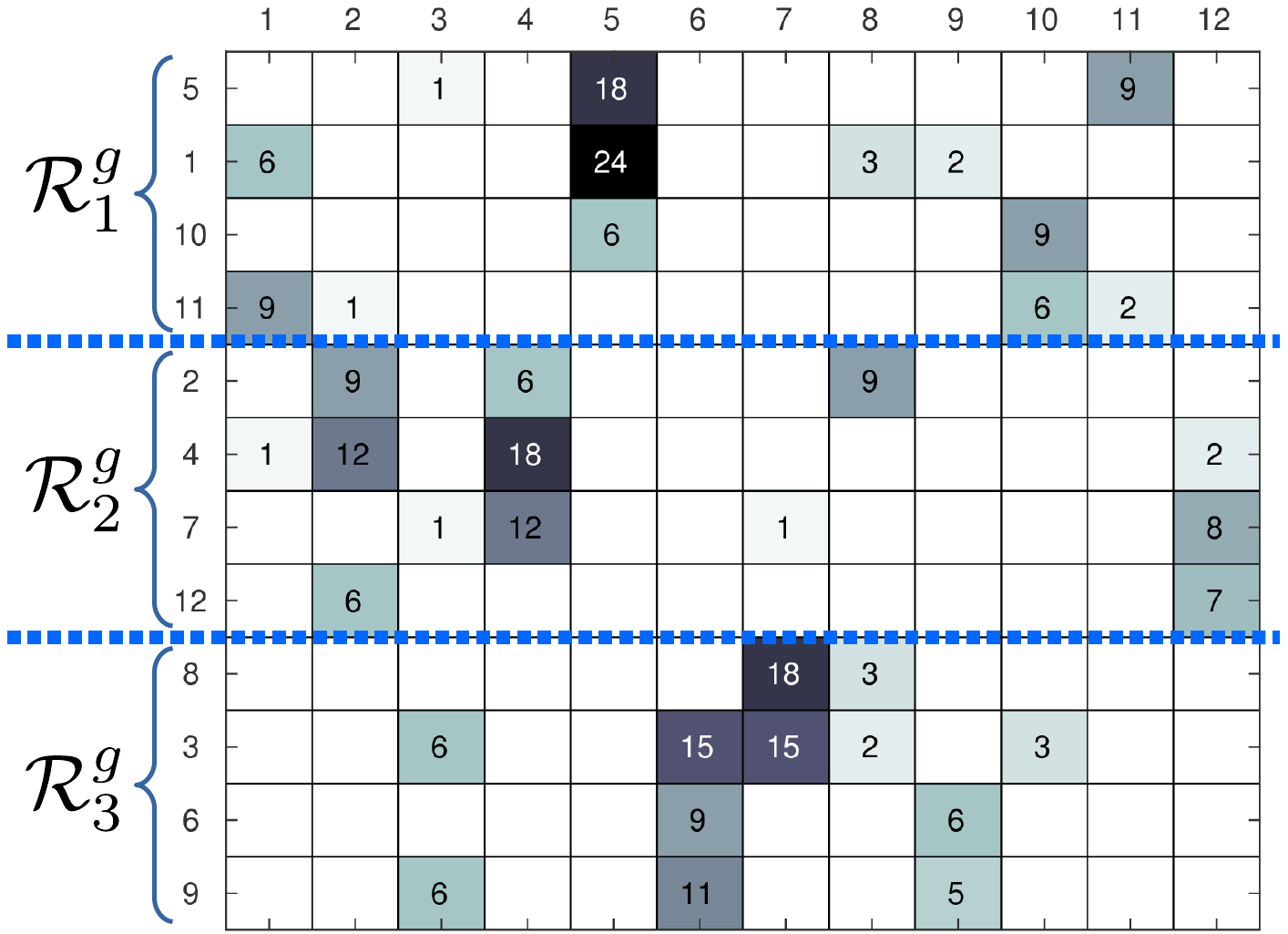}
        \label{metpart}
    }    
 \end{center}  
	\caption{(a)``Good" 3-way  partition $\Pi^g(\mathcal{V})$ of $\mathcal{G}_{\rm{RIP}}(A)$ and (b) 3-way row partition of $A$ induced by~\cref{met}.}
\label{graphpart2}
\vspace{-0.3cm}
\end{figure}

We introduce~\cref{graphpart1,graphpart2} in order to clarify the proposed graph partitioning method for block-row partitioning.
\Cref{unipart} shows a straightforward 3-way block-row partition $  \{ \mathcal{R}_1^{s},\mathcal{R}_2^s,R_3^s \}$ of the sample matrix $A$ given in~\cref{graph11}, where the first four, the second four and the third four consecutive rows in the original order constitute the block rows $\mathcal{R}_1^{s},\mathcal{R}_2^s$ and $\mathcal{R}_3^s$, respectively.
\Cref{uni} shows the $3$-way vertex partition  $\Pi^s(\mathcal{V}) = \{\mathcal{V}_1^s,\mathcal{V}_2^s,\mathcal{V}_3^s\}$ of $\mathcal{G}_{\rm{RIP}}$ that corresponds to this straightforward $3$-way block-row partition.  
\Cref{met} shows a ``good"  $3$-way vertex partition $\Pi^g(\mathcal{V}) = \{\mathcal{V}_1^g,\mathcal{V}_2^g,\mathcal{V}_3^g\}$ of $\mathcal{G}_{\rm{RIP}}$ obtained by using the graph partitioning tool METIS~\cite{metis}.
\Cref{metpart} shows the permuted $A^\Pi$ matrix and block-row partition $ \{\mathcal{R}_1^{g},\mathcal{R}_2^{g},R_3^{g}\}$  induced by the 3-way vertex partition $\Pi^g(\mathcal{V}).$ 

As seen in~\cref{graphpart1,graphpart2},  both straightforward and ``good" block-row partitions achieve perfect balance on the row counts of blocks by having exactly four rows per block.
The quality difference between straightforward and ``good" block-row partitions can be easily seen by comparing the $3$-way partitions of $\mathcal{G}_{\rm{RIP}}$ in \cref{uni,met}, respectively.
As seen in~\cref{uni}, eight out of nine thick edges remain on the cut of $\Pi^s(\mathcal{V})$,  whereas all of the nine thick edges remain internal in $\Pi^g(\mathcal{V})$ as seen in~\cref{met}. 

\begin{figure}[hbpt]
\vspace{-0.3cm}
  \begin{center}  
       \subfloat []{
        \includegraphics[scale=.37]{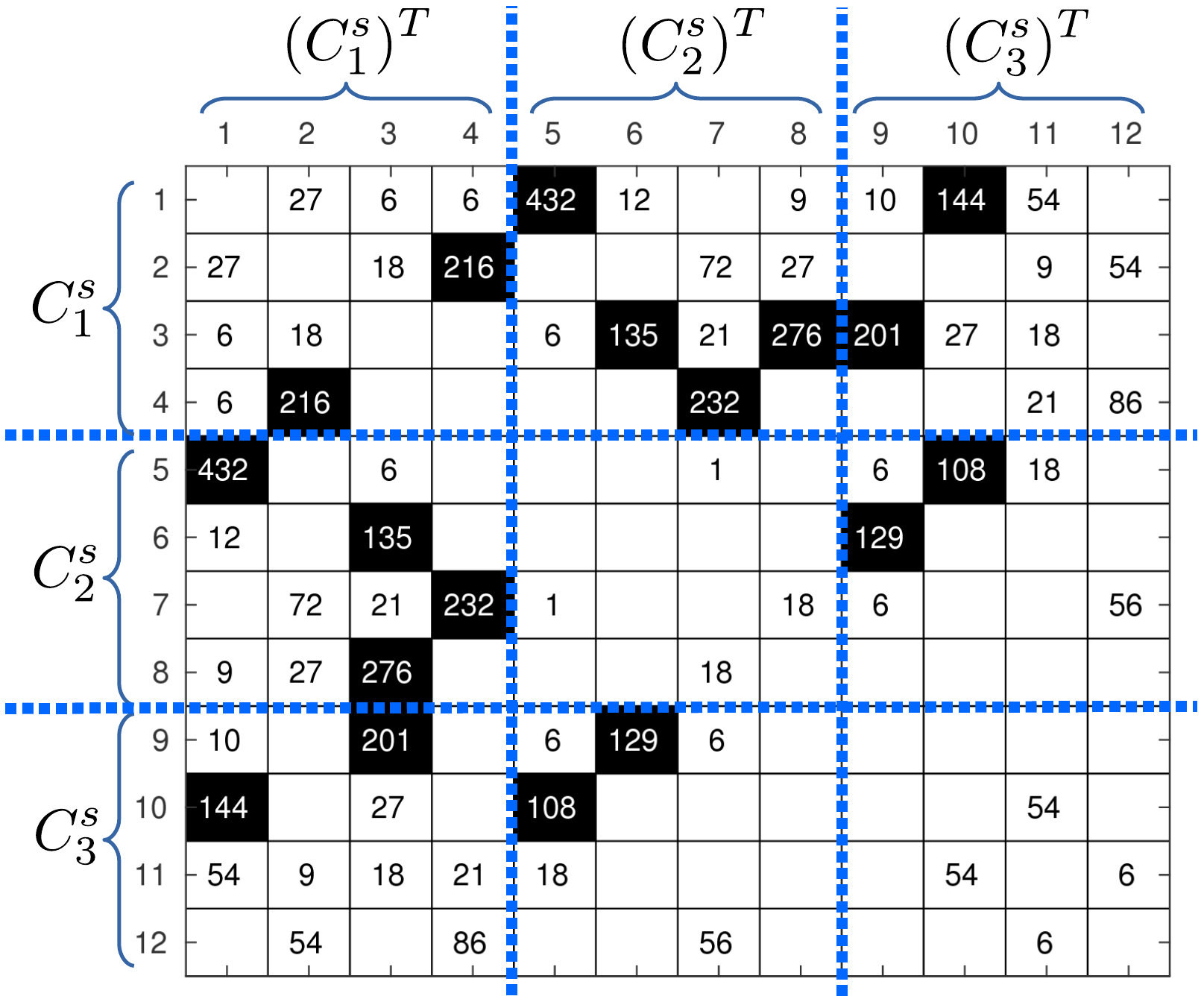}
        \label{aat}
    } 
	\subfloat [] {
        \includegraphics[scale=.37]{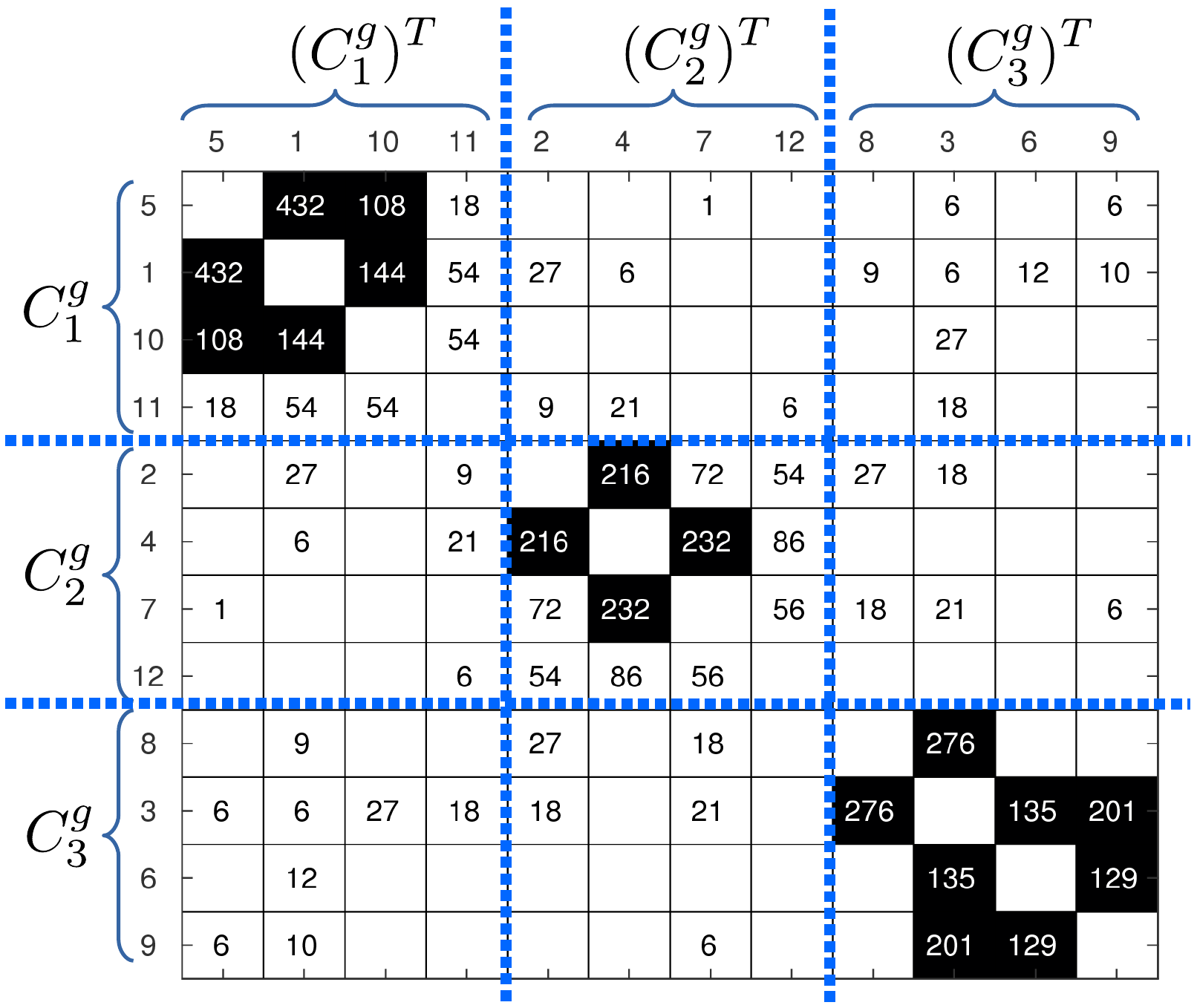}
        \label{aatmetis}
    }  
 \end{center}  
 \caption{$3 \times 3$ block-checkerboard partition of matrix $C=AA^T$  induced by 3-way (a) straightforward (\cref{unipart}) and  (b) ``good" (\cref{metpart}) block-row partitions of matrix $A$.}
\label{sampleAAT}
\vspace{-0.1cm}
\end{figure}

\Cref{sampleAAT} shows  the $3 \times 3$ block-checkerboard partitioning of the resulting matrix $C=AA^T$ induced by straightforward and ``good" block-row partitioning of the sample matrix $A$ in~\cref{unipart,metpart}, respectively.
Note that both rows and columns of the $C$ matrix are partitioned conformably with the row partitions of the $A$ matrix.
The comparison of~\cref{aat,aatmetis} shows that large nonzeros (dark cells) are scattered across the off-diagonal blocks of matrix $C$ for the straightforward partitioning, whereas large nonzeros (dark cells) are clustered to the diagonal blocks of $C$ for the ``good" partitioning.

\begin{figure}[hbpt]
\vspace{-0.3cm}
  \begin{center}  
	  \begin{minipage}{.5\textwidth}
	  \centering  
	  \begin{minipage}{.25\textwidth}
  	  		IP$(\Pi^s)$ =
	  \end{minipage}
	   \begin{minipage}{.5\textwidth}
    \subfloat []{
        \includegraphics[scale=.55]{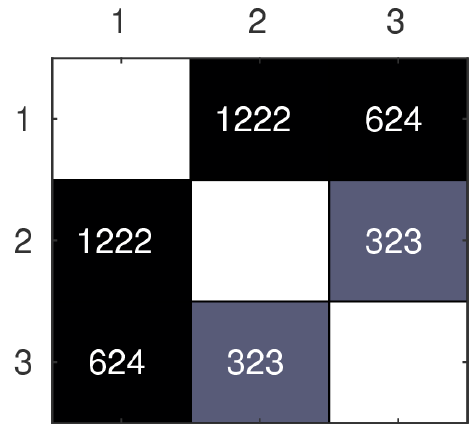}%
        \label{sumprod}
    }%
    \end{minipage}%
    \end{minipage}%
     \begin{minipage}{.5\textwidth}
     	\centering  
       \begin{minipage}{.25\textwidth}
  	  		IP$(\Pi^g)$ =
	  \end{minipage}
	   \begin{minipage}{.5\textwidth}
    \subfloat [] {
        \includegraphics[scale=.55]{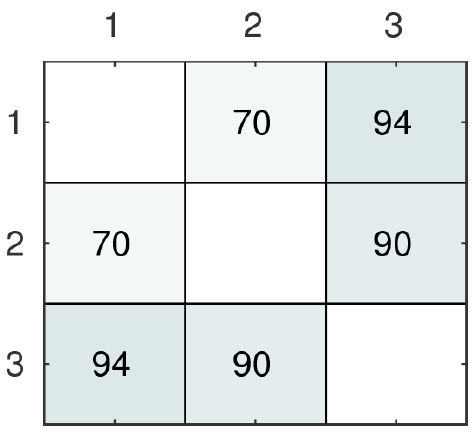}%
        \label{sumprodmetis}
    }	%
    \end{minipage}%
        \end{minipage}
 \end{center}  
 \caption{Inter-block row-inner-product matrix \rm{IP}$(\Pi)$ for (a) straightforward  and (b) ``good" block-row partitions.}
\label{sumprods}
\end{figure}

We introduce~\cref{sumprods} to compare the quality of the straightforward and ``good" block-row partitions in terms of inter-block row-inner-product sums.
In the figure, each off-diagonal entry ip$_{km}$ of the $3 \times 3$ IP matrix shows the sum of the inter-block row inner-products between the respective block rows $\mathcal{R}_k$ and $\mathcal{R}_m$.
That is, 
\begin{equation}
{\mathrm{ip}}_{km} \triangleq \sum\limits_{\substack{ {\mathrm{r}}_i \in \mathcal{R}_k \\ {\mathrm{r}}_j \in \mathcal{R}_m } } \abs{ {\rm{r}}_i{\mathrm{r}}_j^T } \; \text{ for } k \neq m.
\end{equation}
As seen in~\cref{sumprods}, ip$^s_{12} = 1,\!222$ for the straightforward partition, whereas ip$^g_{12}  = 70$ for the ``good" partition.
Note that ip$_{km}$ is also equal to the sum of the absolute values of the nonzeros of the off-diagonal block $C_{km}$ at the $k^{th}$ row block and  $m^{th}$ column block of the $C$ matrix, i.e., 
\begin{equation}
{\rm{ip}}_{km} = \sum\limits_{\substack{{\rm{r}}_i \in R_k \\ {\rm{r}}_j \in R_m }} |c_{ij}| . 
\end{equation}
Therefore the total sum of inter-block inner products is
\begin{equation}
\begin{aligned}
 \mathrm{interIP}(\Pi^s)  &= {\rm{ip}}^s_{12} + {\rm{ip}}^s_{13} + {\rm{ip}}^s_{23} \\
 & =  1,\!222+ 624 + 323 = 2,\!169 \\
\end{aligned}
\end{equation}
for the straightforward partition, whereas for the ``good" partition it is
\begin{equation}
\begin{aligned}
 \mathrm{interIP}(\Pi^g)  = 70 + 94 + 90 = 254 
 \end{aligned}.
\end{equation}

\begin{figure}[tbhp]
	\vspace{-0.4cm}
	\begin{center}    
		\subfloat [Straightforward partitioning (\cref{unipart})] {
			\includegraphics[scale=.55]{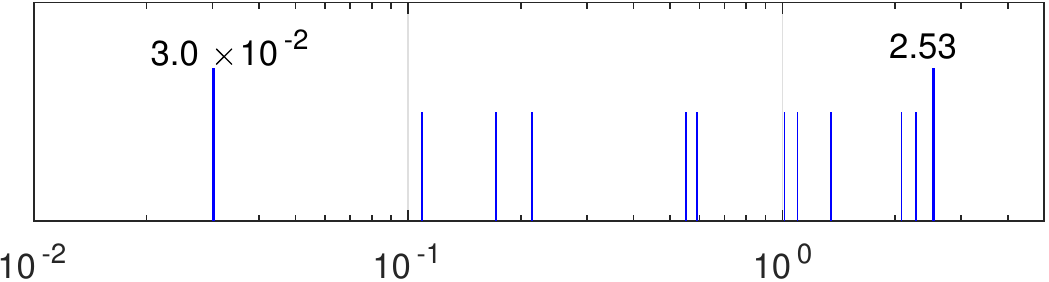}
			\label{eig}
		}	  \hspace{0.2cm}
		\subfloat [``Good" partitioning (\cref{metpart})]{
			\includegraphics[scale=.55]{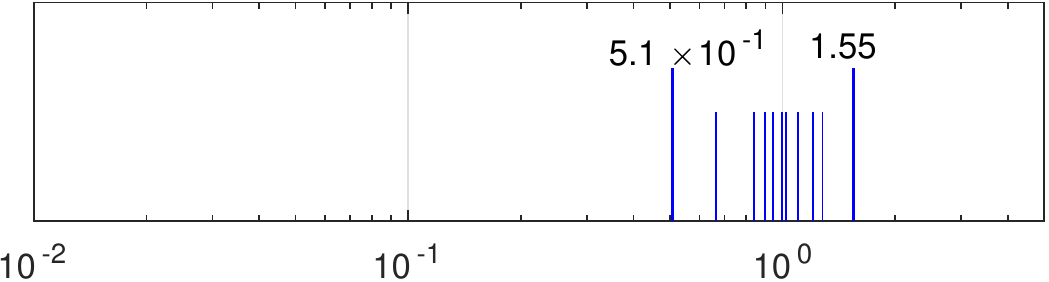}
			\label{metiseig}
		}   
	\end{center}  
	\caption{Eigenvalue spectrum of $H$ for the block-row partitionings of the sample matrix given in~\cref{graph11}.}
	\label{sampleeig}
	\vspace{-0.4cm}
\end{figure}

\Cref{eig,metiseig}, respectively show  the eigenvalue spectrum of $H$ for the straightforward and ``good" partitionings. 
As seen in the figures, for the straightforward partitioning the eigenvalues reside in the interval $[3.0\!\times\!10^{-2},2.53]$, whereas for "good" partitioning the eigenvalues reside in the interval $[5.1\!\times\!10^{-1},1.55]$.
As seen in~\cref{metiseig}, after using $\mathcal{G}_{\mathrm{RIP}}$ partitioning the eigenvalues are much better clustered around $1$ and the smallest eigenvalue is much larger than that of the straightforward partitioning method.

\subsection{Implementation}
\label{subsec:impl}

Implementation of the proposed partitioning method consists of two stages which are constructing $\mathcal{G}_{\rm{RIP}}$ and partitioning $\mathcal{G}_{\rm{RIP}}$.\\ 

\noindent
\textbf{Constructing $\mathcal{G}_{\rm{RIP}}$:} 
For constructing $\mathcal{G}_{\rm{RIP}}$, we choose to use basic sparse matrix-matrix multiplication (SpGEMM)~\cite{gustavson1978two} kernel due to existing efficient implementations.
The edges of the $\mathcal{G}_{\rm{RIP}}$ are obtained from the nonzeros of the $C=AA^T$ matrix, whereas their weights are obtained from the absolute values of those nonzeros.

\begin{figure}[htbp]
\vspace{-0.4cm}
  \begin{center}    
    \subfloat [$A$] {
       \includegraphics[scale=.23]{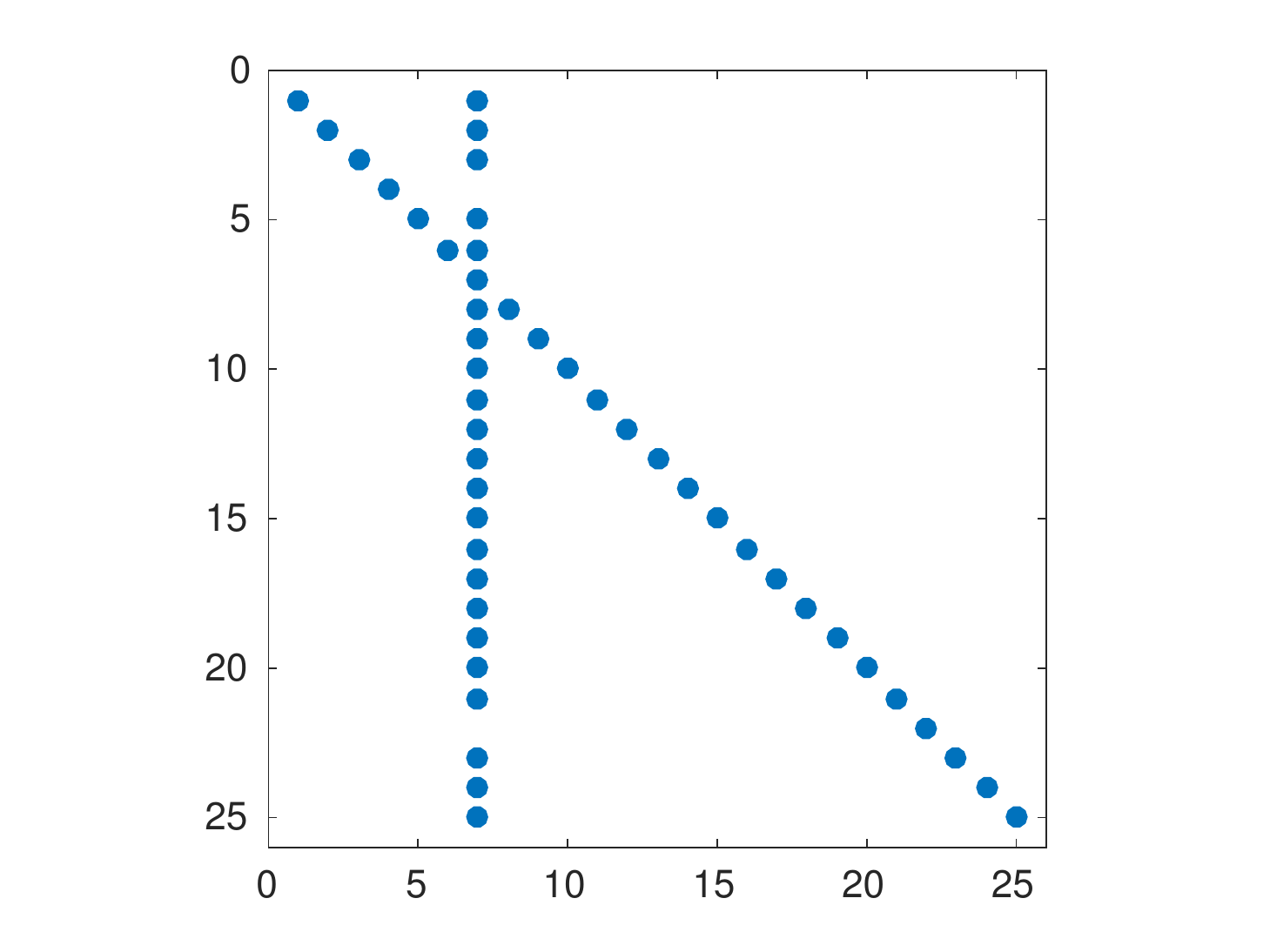}
        \label{eye1}
    }	\hspace{-0.9cm}  
	\subfloat [$C=AA^T$]{
        \includegraphics[scale=.23]{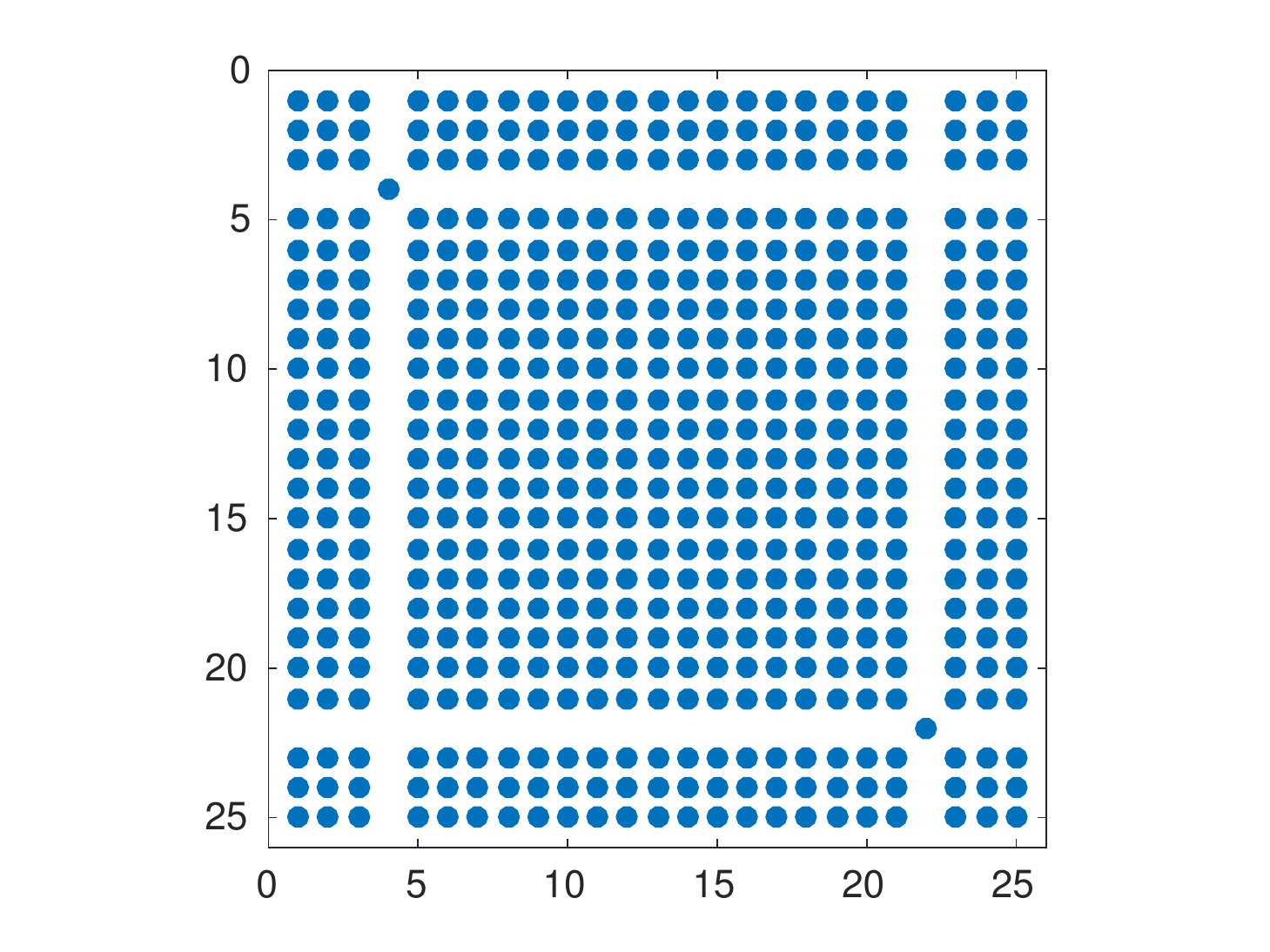}
        \label{eye2}
    }	\hspace{-0.1cm}  
    \subfloat [$\tilde{A}$] {
       \includegraphics[scale=.23]{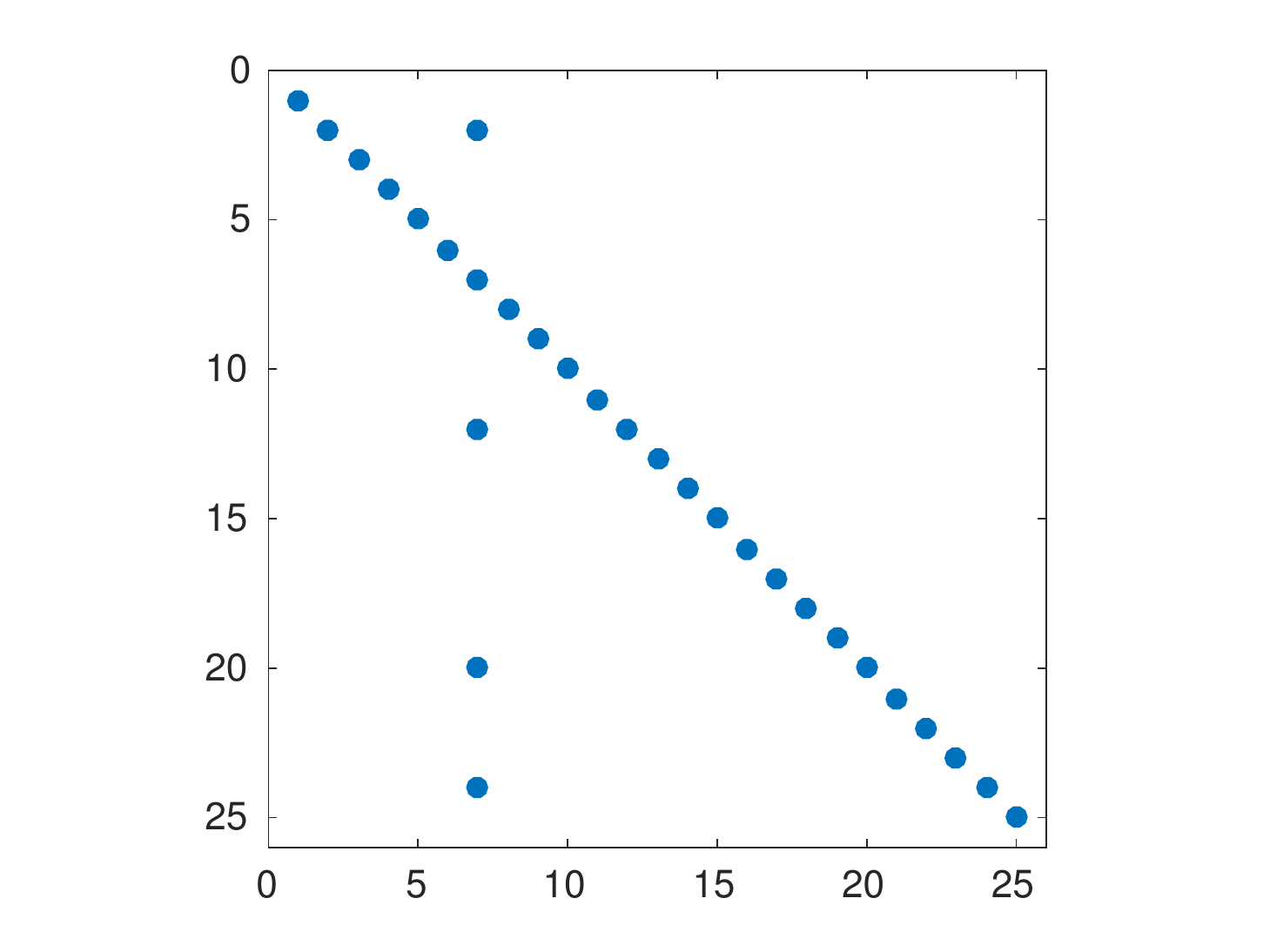}
        \label{eye3}
    }	\hspace{-0.9cm}  
	\subfloat [$\tilde{C}=\tilde{A}\tilde{A}^T$]{
        \includegraphics[scale=.23]{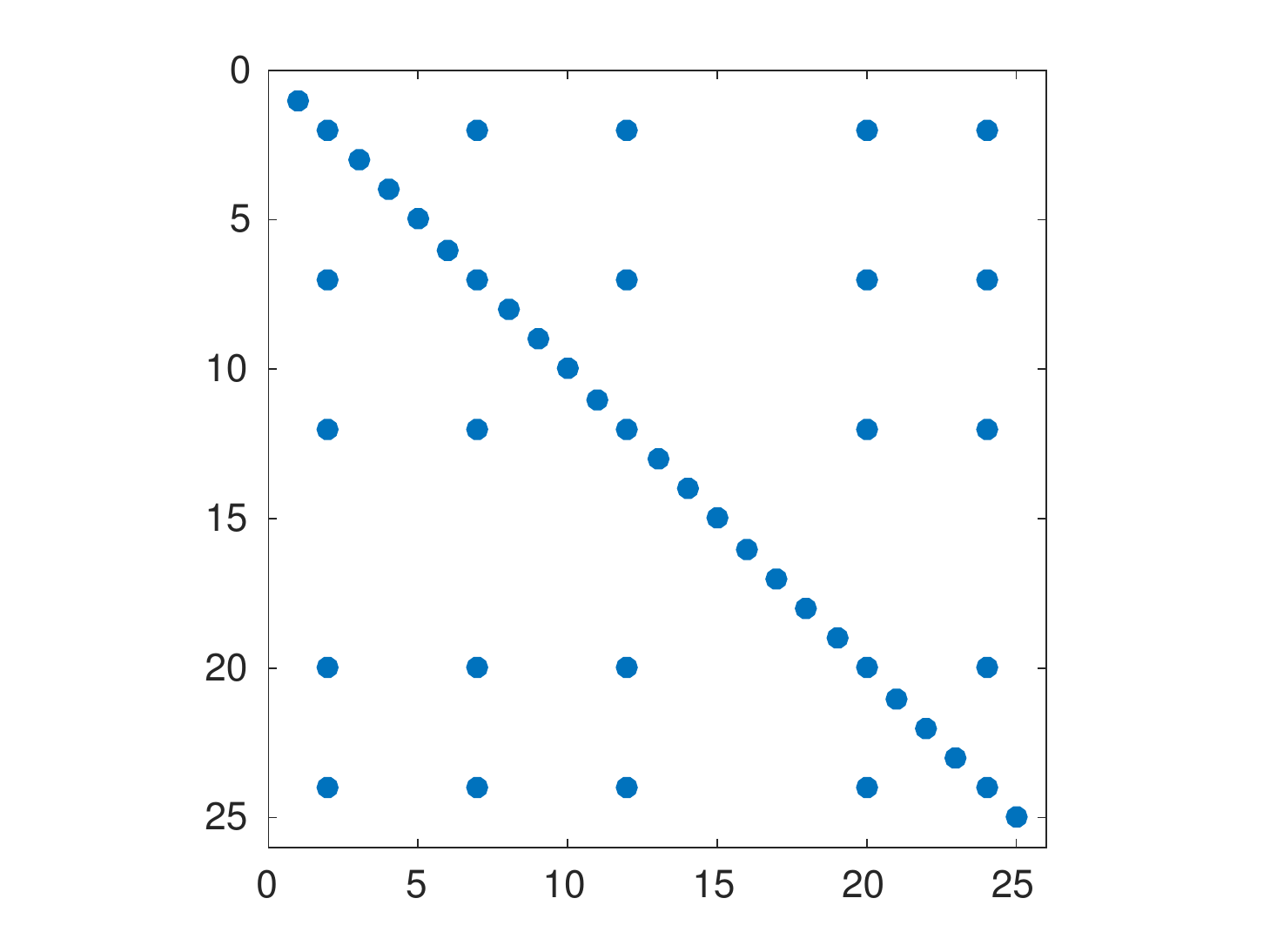}
        \label{eye4}
    }	
 \end{center}  
	\caption{Nonzero patterns of $A$, $AA^T$, $\tilde{A}$ and $\tilde{A}\tilde{A}^T$.}
\label{eye}
\vspace{-0.25cm}
\end{figure}

Note that when matrix $A$ has dense column(s), the corresponding matrix $C =AA^T$ will be quite dense. 
In other words, when a column has $nz$ nonzeros, corresponding $C$ matrix will have at least $nz^2$ nonzeros. 
For example,~\cref{eye1} shows a $25 \times 25$ sparse matrix $A$ which has a dense column having $23$ nonzero entries.
As seen in~\cref{eye2}, matrix $AA^T$ is dense as it has $531$ nonzero entries.
Clearly, large number of nonzeros in $C$ (i.e., large number of edges in $\mathcal{G}_{\rm{RIP}}$) increases the memory requirement and computation cost of SpGEMM as well as the time requirement for partitioning $\mathcal{G}_{\rm{RIP}}$.

In order to alleviate the aforementioned problem, we propose the following met-hodology for sparsifying $C$.
We identify a column $A(:,i)$ (in MATLAB~\cite{MATLAB:2015} notation) of an $n\!\times\!n$ matrix $A$ as a dense column if it contains more than $\sqrt{n}$ nonzeros \mbox{($nnz(A(:,i)) > \sqrt{n}$)}.
Given $A$, we extract a sparse matrix $\tilde{A}$ by keeping the largest (in absolute value) $\sqrt{n}$ nonzeros of dense columns of $A$.
That is,  the smallest $nnz(A(:,i)) - \sqrt{n}$  entries of a dense $A$-matrix column $A(:,i)$ is ignored during constructing column $\tilde{A}(:,i)$ of $\tilde{A}$.
Hence, the SpGEMM operation is performed on $\tilde{A}$ to obtain sparsified resulting matrix $\tilde{C}=\tilde{A}\tilde{A}^T$.
This will lead to a sparsified $\mathcal{\tilde{G}}_{\rm{RIP}}$ graph.
For example,~\cref{eye3} shows this sparsity pattern of a sparse matrix $\tilde{A}$ which is extracted from $A$ by keeping $5$ largest nonzeros in the dense column of $A$.
As seen in~\cref{eye4}, matrix $\tilde{C}=\tilde{A}\tilde{A}^T$ is very sparse with respect to~\cref{eye2}.
Note that both $\tilde{A}$ and $\tilde{C}$ are used only for constructing  $\mathcal{\tilde{G}}_{\rm{RIP}}$ of $A$. After the partitioning stage, both matrices are discarded.
In the rest of the paper, $\mathcal{\tilde{G}}_{\rm{RIP}}$ will be referred to as $\mathcal{G}_{\rm{RIP}}$ for the sake of the simplicity of presentation.

\noindent \textbf{{Partitioning $\mathcal{G}_{\rm{RIP}}$}:}
We use multilevel graph partitioning tool METIS~\cite{metis} for partitioning $\mathcal{G}_{\rm{RIP}}$.
In order to compute integer edge weights required by METIS, we multiply the floating-point edge cost values with $\alpha$ and round them up to the nearest integer value; $wgt(v_i,v_j)  =  \ceil*{ \alpha \times {cost(v_i,v_j) }},$  where $\alpha$ is  a sufficiently large integer.
Here, $cost(v_i,v_j)$ is the edge cost computed according to~\cref{eq14} and $wgt(v_i,v_j)$ is the weight of the respective edge provided to METIS.
Since the rows of matrix $A$ are prescaled to have 2-norm equal to one in the preprocessing phase, each edge cost $cost(v_i,v_j)$ should be in the range $(0,1]$ and the resulting edge weight $wgt(v_i,v_j)$ will be an integer in the range $[1,\alpha]$.


\section{Experimental Results}
\label{section3}

\subsection{Experimental Framework}
In the experiments, we used the CG accelerated block Cimmino implementation available in ABCD Solver v1.0~\cite{abcd}.
In ABCD Solver, we used MUMPS 5.1.2~\cite{MUMPS} sparse direct solver to factorize the systems in~\cref{augmentedMatrix} once and solve the system iteratively.
We note that the proposed scheme is designed for the classical block Cimmino by improving the numerical orthogonality between blocks and it does not intend to improve the structural orthogonality. 
Hence, it is not applicable to the augmented block Cimmino algorithm that is also available in the ABCD Solver where the number of the augmented columns depends only on structural orthogonality.

We adopted the same matrix scaling procedure~\cite{amestoy2008parallel} as in ABCD Solver.
This is a parallel iterative procedure which scales the columns and rows of $A$ so that the absolute value of largest entry in each column and row is one. 
We first perform row and column scaling in order to avoid problems due to poor scaling of the input matrix.
Then, we also perform row scaling on $A$ to have 2-norm equal to exactly one, so that the actual values in $AA^T$ would then correspond to cosines of the angles between pairs of rows in $A$~\cite{drummond2015partitioning,zenadi2013methodes}. 
We note that  $H$  is numerically independent of row scaling. However the column scaling affects $H$ and it can be considered as  preconditioner~\cite{ruiz1992solution,zenadi2013methodes}.

ABCD Solver includes a stabilized block-CG accelerated block Cimmino algorithm~\cite{arioli1995block} especially for solving systems with multiple right-hand side vectors.
Since the classical CG is guaranteed to converge~\cite{o1980block} for systems where the coefficient matrix is symmetric and positive definite and its convergence theory is well-established,  in this work we utilized the classical CG accelerated block Cimmino in~\cref{algo2} for solving sparse linear system of equations with single right-hand side vector rather than the block-CG acceleration. 

In parallel CG accelerated block Cimmino algorithm, the work distribution among processors is performed in exactly the same way as in ABCD Solver.
That is, if the number of row blocks is larger than the number of processors, row-blocks are distributed among processors so that each processor has equal workload in terms of number of rows.
If the number of row-blocks is smaller than the number of processors, master-slave computational approach~\cite{arioli1995parallel,drummond2015partitioning,duff2015augmented} is adopted.
Each master processor owns a distinct row-block and responsible for inner-product and matrix-vector computations.
Each slave processor is a supplementary processor which helps the specific master processor in the factorization and solution steps of MUMPS.
After the analysis phase of MUMPS, slave processors are mapped to some master processors according to the information of FLOP estimation in the analysis phase.

In all experiments with the CG accelerated block Cimmino, we use the normwise backward error~\cite{arioli1992stopping} at iteration $t$
\begin{equation}
\label{stopping}
 \gamma^{(t)} = \dfrac{\Vert Ax^{(t)} -f \Vert_\infty}{ \Vert A \Vert_\infty \Vert x^{(t)} \Vert_1 + \Vert f \Vert_\infty}  < 10^{-10}
\end{equation}
as the stopping criterion and we use $10,\!000$ as the maximum number of iterations.
The right-hand side vectors of the systems are obtained by multiplying the coefficient matrices with a vector whose elements are all one.
In all instances, the CG iterations are started from the zero vector.

In the experiments, we have compared the performance of the proposed partitioning method (GP) against two baseline partitioning methods already available in ABCD Solver.
The first baseline method, which is referred to as uniform  partitioning (UP) method in ABCD Solver, partitions the rows of the coefficient matrix into a given number of block rows with almost equal number of rows without any permutation on the rows of $A$.
Note that UP method is the same as the straightforward partitioning method mentioned in \cref{section2}.

The second baseline method is the hypergraph partitioning (HP) method~\cite{drummond2015partitioning}.
This method uses the column-net model~\cite{1999hypergraph} of sparse matrix $A$ in which rows and columns are respectively represented by vertices and hyperedges both with unit weights~\cite{drummond2015partitioning}.
Each hyperedge connects the set of vertices corresponding to the rows that have a nonzero on the respective column.
In the HP method, a $K$-way partition of the vertices of the column-net model is used to find  $K$ block rows.
The partitioning constraint of maintaining balance on part weights corresponds to finding block rows with equal number of rows.

CM-based partitioning strategies~\cite{drummond2015partitioning,ruiz1992solution,zenadi2013methodes} are also considered as another baseline approach.
However, our experiments showed that CM-based strategies fail in producing the desired number of balanced partitions and achieving convergence for most of the test instances.
Due to a significantly larger number of failures of CM-based strategies compared to UP and HP, UP and HP are selected as a much better baseline algorithms.

The HP method in ABCD Solver uses the multilevel hypergraph partitioning tool PaToH~\cite{1999patoh} for partitioning the column-net model of matrix $A$. 
Here, we use the same parameters for PaToH specified in ABCD Solver. 
That is, the final imbalance (i.e., $\epsilon$ in \cref{partbal}) and initial imbalance (imbalance ratio of the coarsest hypergraph) parameters in PaToH are set to $50\%$ and $100\%$, respectively.
The other parameters are left as default of PaToH as in ABCD Solver.

Parallel block Cimmino solution times on some real problems are experimented by using PaToH with different imbalance ratios in \cite{drummond2015partitioning}. 
It is reported that although partitioning with weak balancing reduces greatly the number of interconnections which leads to decrease in the number of iterations, however, it increases the parallel solution time because of highly unbalanced computational workload among processors. 
Therefore, finding ``good" partition imbalance ratios can be important for the parallel performance of block Cimmino. 
Due to the space limitation, the impact of different imbalance ratios on the parallel performance is left as a future work.

In the proposed GP method, as mentioned in \cref{subsec:impl}, the multilevel graph partitioning tool METIS is used to partition the row inner-product graph model $\mathcal{G}_{\rm{RIP}}$ of $A$.
The imbalance parameter of METIS is set to 10\% and k-way option is used.
For the sake of a fair comparison between HP and GP methods, unit vertex weights are used in $\mathcal{G}_{\rm{RIP}}$.
The other parameters are left as default of METIS.
Since both PaToH and METIS use randomized algorithms, we report the results of the geometric mean of $5$ experiments with different seeds for each instance.

Here and hereafter, we use GP, HP and UP to refer to the respective block-row partitioning method as well as ABCD Solver that utilizes the regular block Cimmino algorithm for solving the systems partitioned by the respective method. 

The extensive numerical experiments were conducted on a shared memory system. 
The shared memory system is a four socket 64-core computer that contains four AMD Opteron 6376 processors, where each  processor has 16 cores running at 2.3GHz and a total of 64GB of DDR3 memory. 
Due to memory bandwidth limitations of the platform, experiments are performed with 32 cores. 

ABCD Solver~\cite{abcd} is implemented in C/C++ programming language with MPI-OpenMP based hybrid parallelism. 
Furthermore, an additional parallelism level can be incorporated with multithreaded BLAS/LAPACK libraries.
However, in the experiments, we used pure MPI-based parallelism which gives the best performance for our computing system.

\subsection{Effect of block-row partitioning on eigenvalue spectrum of $H$}
The convergence rate of the CG accelerated block Cimmino algorithm is related to the eigenvalue spectrum of $H$.   
By the nature of the block Cimmino algorithm, most of the eigenvalues of  $H$  are clustered around 1, but there can be some eigenvalues at extremes of the spectrum. 

In this subsection, we conduct experiments to study the effect of the partitioning on the eigenvalue spectrum of $H$ by comparing the proposed GP method against UP and HP. 
In these experiments,  in order to be able to compute the eigenvalue spectrum requiring reasonable amount of time and memory, we use four small nonsymmetric sparse matrices: ${\mathtt{sherman3}}$, ${\mathtt{GT01R}}$,    ${\mathtt{gemat11}}$ and ${\mathtt{LeGresley\_4908}}$ 
from SuiteSparse Matrix Collection~\cite{davis2011university}.
The first and the second matrices arise in Computational Fluid Dynamics problems, whereas the third and fourth matrices arise in Power Network problems.
We partition the matrices into $8$ block rows for all of the three partitioning methods UP, HP and GP.

\begin{figure}[htbp]
	{
		\tiny
		\begin{center}    
			\subfloat [$\mathtt{sherman3}:$ $n=5,\!005$, $nnz=20,\!033$] {
				\includegraphics[scale=.56,trim={1.7cm 0 1.1cm 0.22cm}, clip]{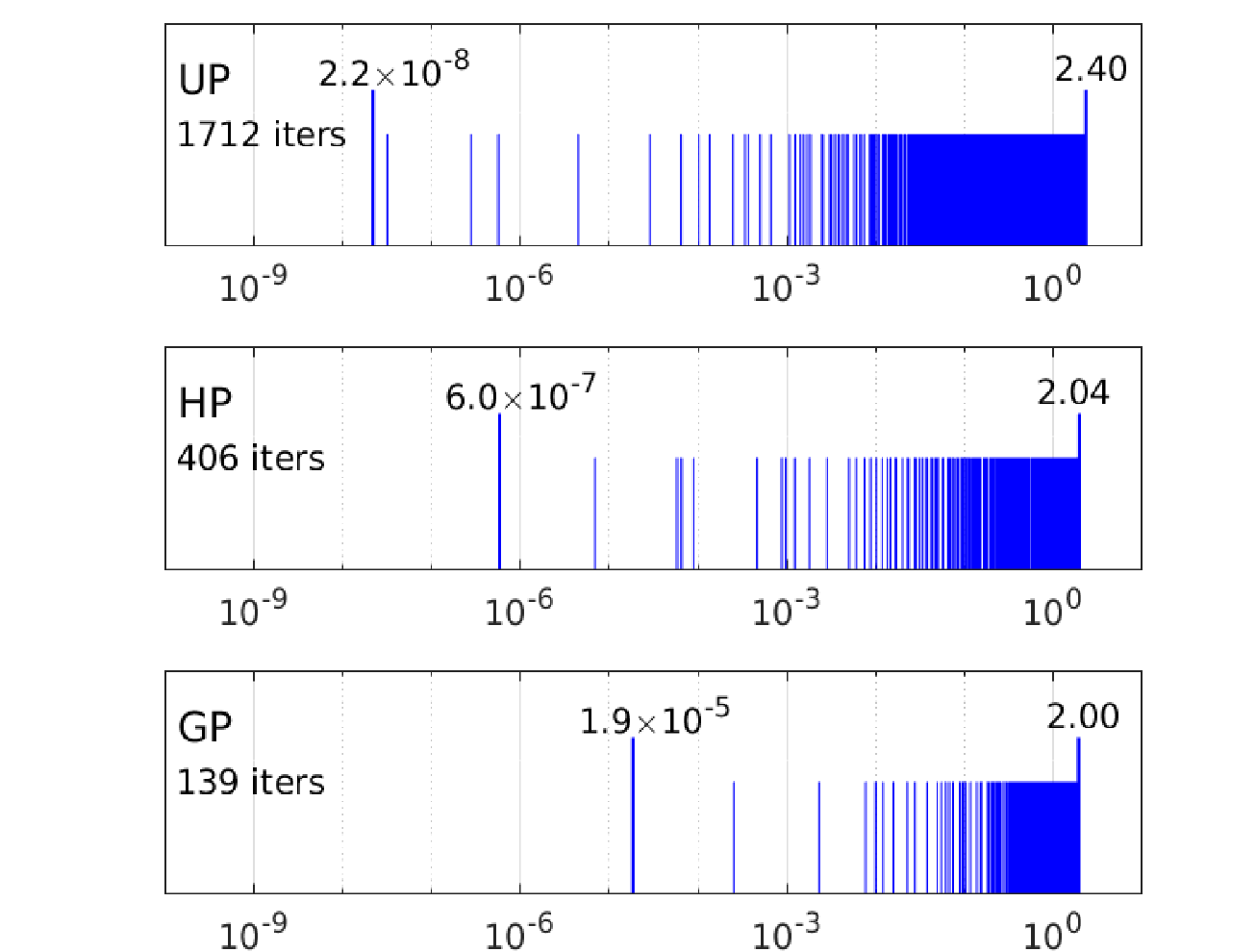}
				\label{sherman_1}
			} \hspace{0.7cm}
			\subfloat [$\mathtt{GT01R}:$ $n=7,\!980$, $nnz=430,\!909$ ] {
				\includegraphics[scale=.56,trim={1.7cm 0 1.1cm 0.22cm}, clip]{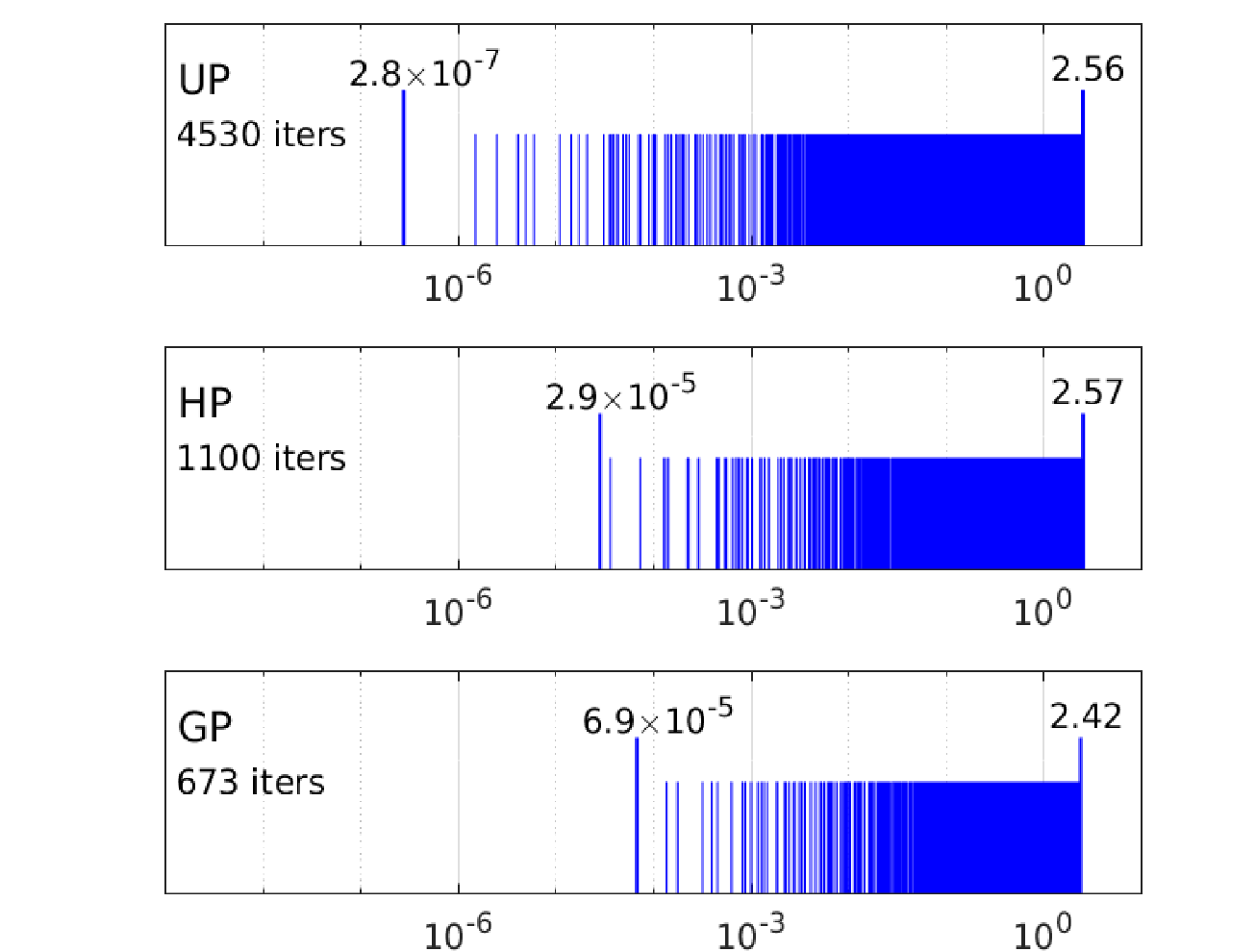}
				\label{GT01R_1}
			}	\\    [+2ex]
			\subfloat [$\mathtt{gemat11}:$ $n=4,\!929$, $nnz=33,\!108$] {
				\includegraphics[scale=.56,trim={1.7cm 0 1.1cm 0.22cm}, clip]{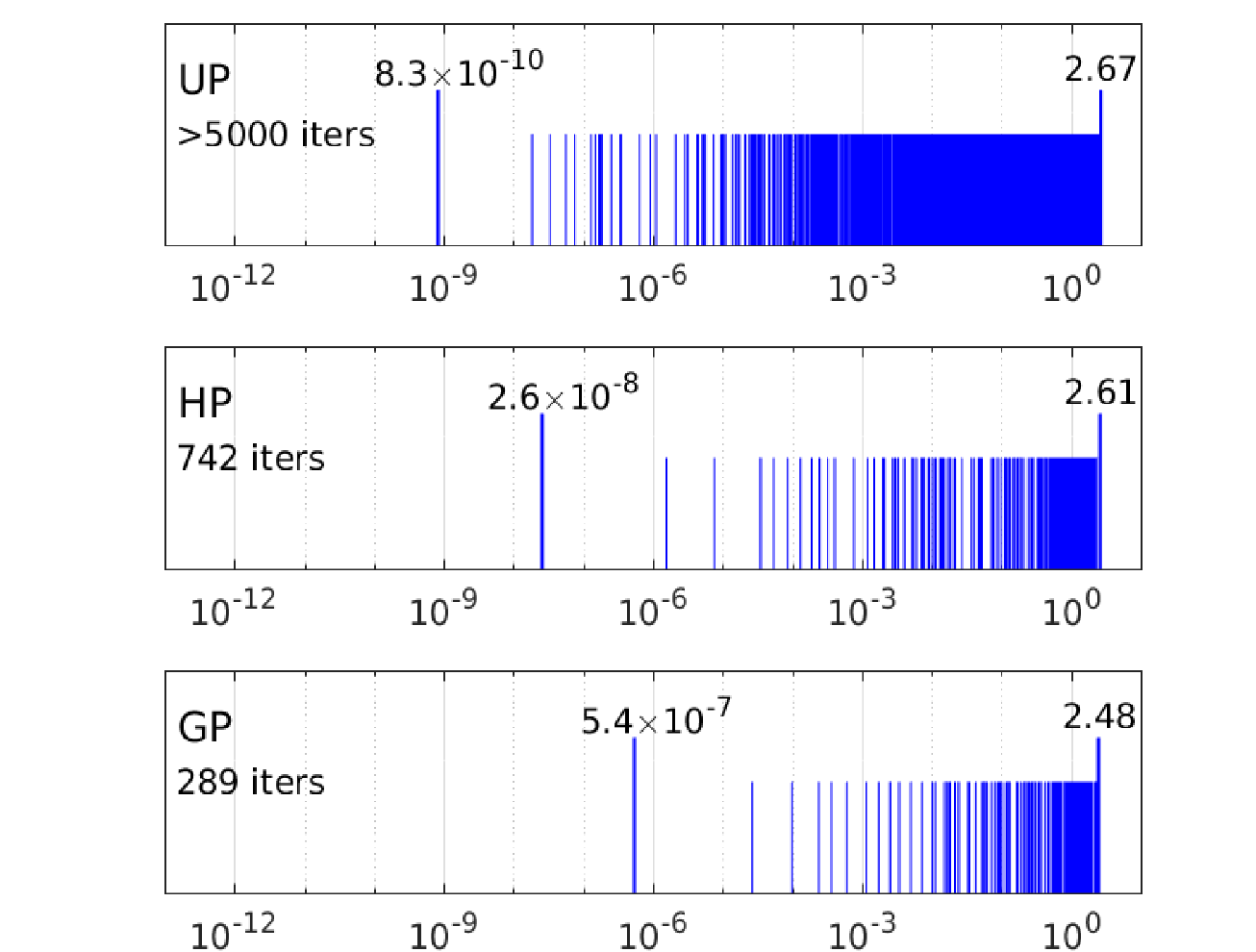}
				\label{LeGresley_1}
			} 	\hspace{0.7cm}     
			\subfloat [$\mathtt{LeGresley\_4908}:$ $n=4\!,\!908$, $nnz=30\!,\!482$] {
				\includegraphics[scale=.56,trim={1.7cm 0 1.1cm 0.22cm}, clip]{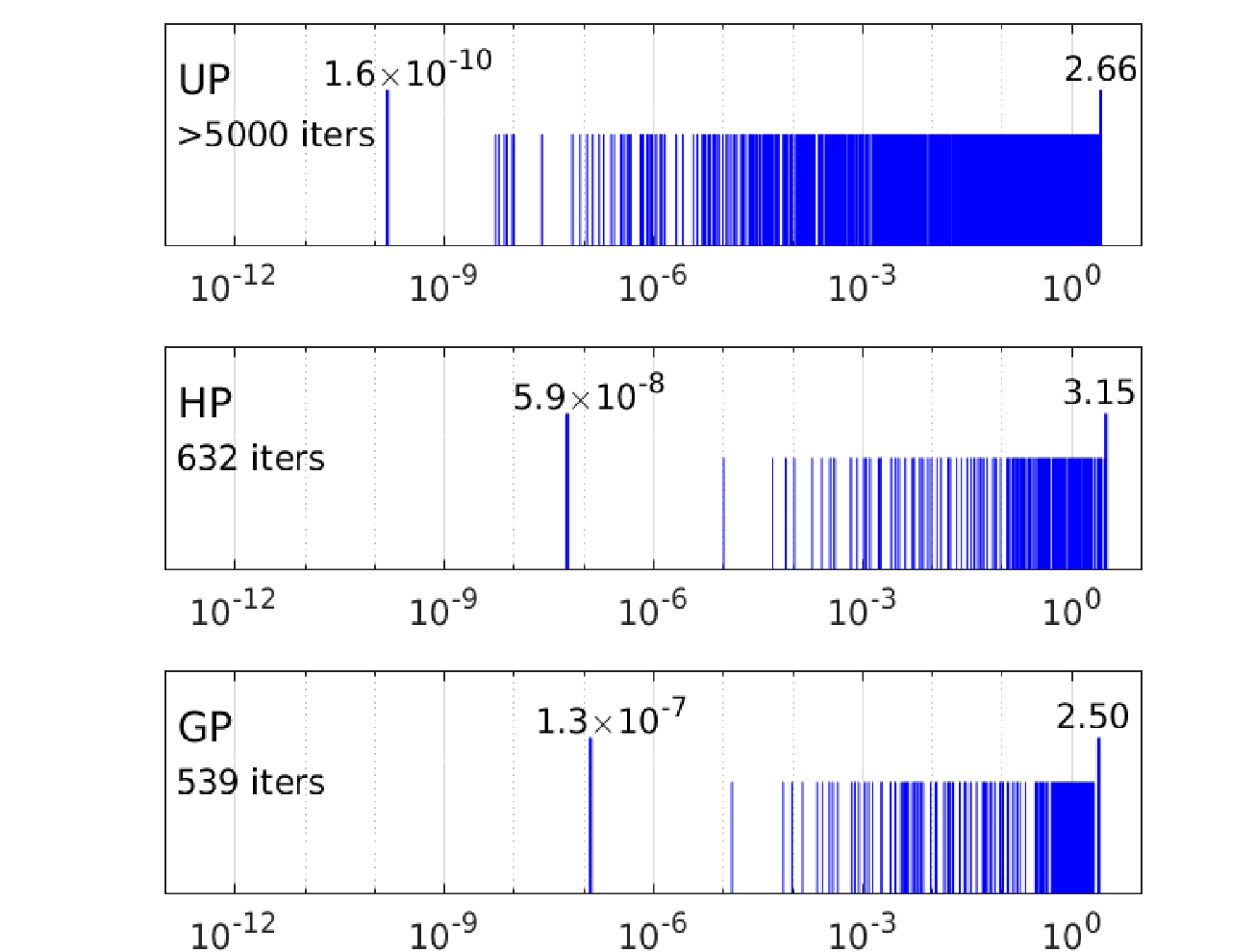}
				\label{gemat11_1}
			}	 	
		\end{center}  
	}
	\caption{Eigenvalue spectrum of $H$ (with the smallest and largest eigenvalues) and the number of CG iterations (iters) required for convergence.}
	\label{eigens}
	\vspace{-0.4cm}
\end{figure}

\Cref{eigens} shows the eigenvalue spectrum of $H$  obtained by UP, HP and GP methods for each test instance.
The figure also reports the number of CG iterations (iters) required for convergence as well as the smallest and largest eigenvalues of $H$.
As seen in the figure, both HP and GP methods achieve significantly better clustering of the eigenvalues around 1 compared to UP.
This experimental finding is reflected to remarkable decrease in the number of CG iterations attained by both HP and GP over UP.

In the comparison of HP and GP, GP achieves better eigenvalue clustering and hence better convergence rate than HP for all instances.
For ${\mathtt{sherman3}}$, ${\mathtt{GT01R}}$, ${\mathtt{gemat11}}$ and ${\mathtt{LeGresley\_4908}}$ instances,
the better clustering quality attained by GP over HP leads to significant improvement in the convergence rate by $66\%$, $39\%$, $61\%$  and $15\%$, respectively.

\subsection{Dataset for Performance Analysis}

For the following experiments, we selected all nonsingular nonsymmetric square matrices whose dimensions are between $50,\!000$ and $5,\!000,\!000$ rows and columns from the SuiteSparse Matrix Collection~\cite{davis2011university}.
The number of matrices based on this criteria turns out to be  $112$. 
Only the  {\tt HV15R} and {\tt cage14} matrices are excluded due to memory limitations.
We have observed that at least one of the three partitioning methods was able to converge in $76$ instances out of these $110$ in less than $10,\!000$ CG iterations.
\Cref{matrix_prop} shows the properties of those $76$ matrices that are used in the experiments.
We note that the main advantage of the block Cimmino algorithm is its amenability to parallelism and requirement of less storage compared to direct methods. It would not be competitive for the smallest problems in the dataset against direct solvers or classical preconditioned iterative solvers.

In \cref{matrix_prop}, the matrices are displayed in increasing sorted order according to their sizes.
The matrices are partitioned into a number of partitions where each row block has approximately $20,\!000$ rows if $n\!>\!100,\!000$ or  $10,\!000$ rows if $n\!<\!100,\!000$.
Thus in our dataset, the smallest matrix is partitioned into 6 row blocks whereas the largest matrix is partitioned into 235 row blocks.

\begin{table}[tbhp]
	\centering
	\caption{Matrix properties ($n$: number of rows/columns, $nnz$:  number of nonzeros)}
	\label{matrix_prop}
	\small
	\resizebox{\textwidth}{!}
	{
		\begin{tabular}{lrr|lrr}
			\hline \vspace{-0.3cm}\\
			\multicolumn{1}{l}{ Matrix name}        & \multicolumn{1}{c}{$n$} & \multicolumn{1}{c}{$nnz$}  & \multicolumn{1}{l}{ Matrix name}        & \multicolumn{1}{c}{$n$} &\multicolumn{1}{c}{$nnz$} \\  
			\hline
			\vspace{-0.3cm}\\
			rajat26             & 51,032  & 247,528   & dc3               & 116,835   & 766,396    \\
			ecl32               & 51,993  & 380,415   & trans4            & 116,835   & 749,800    \\
			2D\_54019\_highK    & 54,019  & 486,129   & trans5            & 116,835   & 749,800    \\
			bayer01             & 57,735  & 275,094   & matrix-new\_3     & 125,329   & 893,984    \\
			TSOPF\_RS\_b39\_c30 & 60,098  & 1,079,986 & cage12            & 130,228   & 2,032,536  \\
			venkat01            & 62,424  & 1,717,792 & FEM\_3D\_thermal2 & 147,900   & 3,489,300  \\
			venkat25            & 62,424  & 1,717,763 & para-4            & 153,226   & 2,930,882  \\
			venkat50            & 62,424  & 1,717,777 & para-10           & 155,924   & 2,094,873  \\
			laminar\_duct3D     & 67,173  & 3,788,857 & para-5            & 155,924   & 2,094,873  \\
			lhr71c              & 70,304  & 1,528,092 & para-6            & 155,924   & 2,094,873  \\
			shyy161             & 76,480  & 329,762   & para-7            & 155,924   & 2,094,873  \\
			circuit\_4          & 80,209  & 307,604   & para-8            & 155,924   & 2,094,873  \\
			epb3                & 84,617  & 463,625   & para-9            & 155,924   & 2,094,873  \\
			poisson3Db          & 85,623  & 2,374,949 & crashbasis        & 160,000   & 1,750,416  \\
			rajat20             & 86,916  & 604,299   & majorbasis        & 160,000   & 1,750,416  \\
			rajat25             & 87,190  & 606,489   & ohne2             & 181,343   & 6,869,939  \\
			rajat28             & 87,190  & 606,489   & hvdc2             & 189,860   & 1,339,638  \\
			LeGresley\_87936    & 87,936  & 593,276   & shar\_te2-b3      & 200,200   & 800,800    \\
			rajat16             & 94,294  & 476,766   & stomach           & 213,360   & 3,021,648  \\
			ASIC\_100ks         & 99,190  & 578,890   & torso3            & 259,156   & 4,429,042  \\
			ASIC\_100k          & 99,340  & 940,621   & ASIC\_320ks       & 321,671   & 1,316,085  \\
			matrix\_9           & 103,430 & 1,205,518 & ASIC\_320k        & 321,821   & 1,931,828  \\
			hcircuit            & 105,676 & 513,072   & ML\_Laplace       & 377,002   & 27,582,698 \\
			lung2               & 109,460 & 492,564   & RM07R             & 381,689   & 37,464,962 \\
			rajat23             & 110,355 & 555,441   & language          & 399,130   & 1,216,334  \\
			Baumann             & 112,211 & 748,331   & CoupCons3D        & 416,800   & 17,277,420 \\
			barrier2-1          & 113,076 & 2,129,496 & largebasis        & 440,020   & 5,240,084  \\
			barrier2-2          & 113,076 & 2,129,496 & cage13            & 445,315   & 7,479,343  \\
			barrier2-3          & 113,076 & 2,129,496 & rajat30           & 643,994   & 6,175,244  \\
			barrier2-4          & 113,076 & 2,129,496 & ASIC\_680k        & 682,862   & 2,638,997  \\
			barrier2-10         & 115,625 & 2,158,759 & atmosmodd         & 1,270,432 & 8,814,880  \\
			barrier2-11         & 115,625 & 2,158,759 & atmosmodj         & 1,270,432 & 8,814,880  \\
			barrier2-12         & 115,625 & 2,158,759 & Hamrle3           & 1,447,360 & 5,514,242  \\
			barrier2-9          & 115,625 & 2,158,759 & atmosmodl         & 1,489,752 & 10,319,760 \\
			torso2              & 115,967 & 1,033,473 & atmosmodm         & 1,489,752 & 10,319,760 \\
			torso1              & 116,158 & 8,516,500 & memchip           & 2,707,524 & 13,343,948 \\
			dc1                 & 116,835 & 766,396   & circuit5M\_dc     & 3,523,317 & 14,865,409 \\
			dc2                 & 116,835 & 766,396   & rajat31           & 4,690,002 & 20,316,253 \\
			\hline
		\end{tabular}
	}
	\vspace{-0.3cm}
\end{table}

\subsection{Convergence and parallel performance}

In this subsection, we study the performance of the proposed GP method against UP and HP in terms of the number of CG iterations and parallel CG time to solution.

\Cref{itertable} shows the number of CG iterations and parallel CG time for each matrix.
In the table, ``$\mathtt{F}$" denotes that an algorithm fails to reach the desired backward error in $10,\!000$ iterations for the respective matrix instance.
As seen in~\cref{itertable}, UP and HP fail to converge in $26$ and $18$ test instances, respectively, whereas GP does not fail in any test instance.

In~\cref{itertable}, the best result for each test instances is shown in bold.
As seen in the table, out of 76 instances, the proposed GP method achieves the fastest convergence in 58 instances, whereas HP and UP achieve the fastest convergence in only 8 and 11 instances, respectively.
As also seen in~\cref{itertable}, GP achieves the fastest iterative solution time in 56 instances, whereas HP and UP achieve the fastest solution time in 11 and 9 instances, respectively.

\begin{table}[tbhp]
	\centering
	\caption{Number of CG iterations and  parallel CG times in seconds }
	\label{itertable}
	\centering
	\scriptsize\renewcommand{\arraystretch}{.9}
	{
		\begin{tabular}{lrrrrrrr}
			\hline \vspace{-0.15cm}\\
			\multirow{3}{*}{Matrix name}        & \multicolumn{3}{l}{\# of CG iterations} &  & \multicolumn{3}{l}{Parallel CG time to soln.} \\   
			
			\cmidrule(l){2-4}  \cmidrule(l){6-8}  
			& UP   & HP    & GP    &  & UP   & HP   & GP      \\ \hline
			\vspace{-0.15cm}\\
			rajat26             & F             & 3303          & \textbf{245}  &  & F               & 62.4           & \textbf{8.4}     \\
			ecl32               & 5307          & 1253          & \textbf{314}  &  & 246.0           & 36.8           & \textbf{10.1}    \\
			2D\_54019\_highK    & 42            & \textbf{4}    & 9             &  & 0.9             & \textbf{0.1}   & 0.2              \\
			bayer01             & 2408          & 382           & \textbf{131}  &  & 63.8            & 8.7            & \textbf{3.4}     \\
			TSOPF\_RS\_b39\_c30 & 676           & \textbf{262}  & 473           &  & 24.9            & \textbf{4.6}   & 8.0              \\
			venkat01            & 54            & 37            & \textbf{34}   &  & 1.8             & \textbf{0.9}   & 0.9              \\
			venkat25            & 915           &625  &  \textbf{599}           &  & 30.8            & \textbf{14.3}  & 14.4             \\
			venkat50            & 1609          & 975           & \textbf{970}  &  & 56.2            & \textbf{23.5}  & 23.7             \\
			laminar\_duct3D     & 466           & 630           & \textbf{394}  &  & 29.6            & 39.0           & \textbf{23.0}    \\
			lhr71c              & F             & 6164          & \textbf{4166} &  & F               & 193.3          & \textbf{136.5}   \\
			shyy161             & \textbf{13}   & 20            & 15            &  & \textbf{0.4}    & 0.6            & 0.5              \\
			circuit\_4          & F             & 256           & \textbf{183}  &  & F               & 15.5           & \textbf{12.3}    \\
			epb3                & 2583          & 3089          & \textbf{2318} &  & 87.7            & 100.9          & \textbf{77.9}    \\
			poisson3Db          & 4797          & 983           & \textbf{715}  &  & 715.0           & 51.9           & \textbf{41.0}    \\
			rajat20             & 629           & 641           & \textbf{322}  &  & 81.5            & 40.9           & \textbf{34.6}    \\
			rajat25             & 1172          & 937           & \textbf{448}  &  & 121.0           & 62.6           & \textbf{48.9}    \\
			rajat28             & 556           & 369           & \textbf{207}  &  & 86.3            & 22.3           & \textbf{21.2}    \\
			LeGresley\_87936    & F             & 7625          & \textbf{3102} &  & F               & 266.5          & \textbf{122.2}   \\
			rajat16             & 6834          & 1022          & \textbf{180}  &  & 835.0           & 68.5           & \textbf{20.0}    \\
			ASIC\_100ks         & F             & \textbf{23}   & 45            &  & F               & \textbf{1.0}   & 2.1              \\
			ASIC\_100k          & 76            & 324           & \textbf{49}   &  & 22.7            & 93.5           & \textbf{13.4}    \\
			matrix\_9           & \textbf{3944} & F             & 9276          &  & \textbf{272.0}  & F              & 704.0            \\
			hcircuit            & 2061          & 2477          & \textbf{460}  &  & 591.0           & 136.5          & \textbf{32.7}    \\
			lung2               & \textbf{12}   & \textbf{12}   & 13            &  & 0.9             & \textbf{0.8}   & 0.9              \\
			rajat23             & F             & 7770          & \textbf{501}  &  & F               & 465.6          & \textbf{33.4}    \\
			Baumann             & \textbf{732}  & 1551          & 1340          &  & \textbf{47.4}            & 86.5  & 89.8             \\
			barrier2-1          & F             & F             & \textbf{1219} &  & F               & F              & \textbf{145.1}   \\
			barrier2-2          & F             & F             & \textbf{1024} &  & F               & F              & \textbf{99.8}    \\
			barrier2-3          & F             & F             & \textbf{1013} &  & F               & F              & \textbf{106.1}   \\
			barrier2-4          & F             & F             & \textbf{1360} &  & F               & F              & \textbf{116.0}   \\
			barrier2-10         & F             & F             & \textbf{1206} &  & F               & F              & \textbf{135.8}   \\
			barrier2-11         & F             & F             & \textbf{1169} &  & F               & F              & \textbf{125.7}   \\
			barrier2-12         & F             & F             & \textbf{1139} &  & F               & F              & \textbf{116.1}   \\
			barrier2-9          & F             & F             & \textbf{1306} &  & F               & F              & \textbf{136.2}   \\
			torso2              & 16            & \textbf{14}   & 15            &  & 0.7             & \textbf{0.6}   & 0.7              \\
			torso1              & F             & \textbf{4376} & 9200          &  & F               & \textbf{322.1} & 833.7            \\
			dc1                 & 629           & 2059          & \textbf{83}   &  & 255.0           & 745.4          & \textbf{32.9}    \\
			dc2                 & 478           & 1313          & \textbf{68}   &  & 170.3           & 504.0          & \textbf{22.4}    \\
			dc3                 & 2172          & 3329          & \textbf{90}   &  & 793.0           & 1220.1         & \textbf{31.9}    \\
			trans4              & 292           & 1416          & \textbf{23}   &  & 105.9           & 452.1          & \textbf{5.8}     \\
			trans5              & 1006          & 4533          & \textbf{33}   &  & 368.7           & 1693.1         & \textbf{8.5}     \\
			matrix-new\_3       & F             & 9739          & \textbf{6707} &  & F               & 709.7          & \textbf{519.4}   \\
			cage12              & \textbf{9}    & 12            & 10            &  & \textbf{1.7}    & 3.1            & 2.2              \\
			FEM\_3D\_thermal2   & 67            & 54            & \textbf{45}   &  & 4.9             & 4.3            & \textbf{3.8}     \\
			para-4              & 7675          & F             & \textbf{3546} &  & 1600.0          & F              & \textbf{423.4}   \\
			para-10             & F             & F             & \textbf{5565} &  & F               & F              & \textbf{588.6}   \\
			para-5              & F             & F             & \textbf{5054} &  & F               & F              & \textbf{610.5}   \\
			para-6              & F             & F             & \textbf{5019} &  & F               & F              & \textbf{574.1}   \\
			para-7              & F             & F             & \textbf{4576} &  & F               & F              & \textbf{508.6}   \\
			para-8              & F             & F             & \textbf{4973} &  & F               & F              & \textbf{535.7}   \\
			para-9              & F             & F             & \textbf{5667} &  & F               & F              & \textbf{682.0}   \\
			crashbasis          & 68            & \textbf{17}   & 21            &  & 10.0            & \textbf{1.1}   & 1.4              \\
			majorbasis          & 48            & \textbf{16}   & 18            &  & 6.6             & \textbf{1.0}   & 1.2              \\
			ohne2               & \textbf{2623} & F             & 3881          &  & F               & 2103.1         & \textbf{820.8}   \\
			hvdc2               & F             & 5622          & \textbf{3042} &  & F               & 446.7          & \textbf{262.1}   \\
			shar\_te2-b3        & \textbf{23}   & 27            & 26            &  & 9.3             & 8.9            & \textbf{8.8}     \\
			stomach             & \textbf{8}    & 11            & 9             &  & \textbf{1.1}    & 1.4            & 1.3              \\
			torso3              & 22            & 30            & \textbf{11}   &  & 5.2             & 7.0            & \textbf{3.1}     \\
			ASIC\_320ks         & F             & 20            & \textbf{2}    &  & F               & 4.5            & \textbf{0.6}     \\
			ASIC\_320k          & 114           & 37            & \textbf{11}   &  & 56.9            & 40.4           & \textbf{8.0}     \\
			ML\_Laplace         & 8615          & 9136          & \textbf{8438} &  & \textbf{3890.3} & 4220.8         & 4215.4           \\
			RM07R               & F             & F             & \textbf{3944} &  & F               & F              & \textbf{10200.0} \\
			language            & 974           & 661           & \textbf{453}  &  & 335.0           & 196.4          & \textbf{159.0}   \\
			CoupCons3D          & 277           & 132           & \textbf{107}  &  & 100.0           & 57.8           & \textbf{48.3}    \\
			largebasis          & 1155          & 701           & \textbf{348}  &  & 549.1           & 174.8          & \textbf{88.3}    \\
			cage13              & \textbf{10}   & 13            & 12            &  & \textbf{6.7}    & 14.7           & 10.6             \\
			rajat30             & 157           & 200           & \textbf{61}   &  & 229.0           & 274.5          & \textbf{86.3}    \\
			ASIC\_680k          & 10            & 19            & \textbf{2}    &  & 11.3            & 38.2           & \textbf{3.5}     \\
			atmosmodd           & \textbf{744}  & 2055          & 1183          &  & \textbf{869.2}  & 2115.2         & 1279.6           \\
			atmosmodj           & \textbf{787}  & 2235          & 1253          &  & \textbf{994.0}  & 2323.2         & 1293.1           \\
			Hamrle3             & F             & 2394          & \textbf{2010} &  & F               & 2375.8         & \textbf{2012.9}  \\
			atmosmodl           & 1206          & 805           & \textbf{379}  &  & 1800.6          & 1005.3         & \textbf{484.0}   \\
			atmosmodm           & 1164          & 697           & \textbf{202}  &  & 1730.0          & 871.4          & \textbf{257.6}   \\
			memchip             & 3278          & 1002          & \textbf{379}  &  & 7450.3          & 1589.9         & \textbf{862.3}   \\
			circuit5M\_dc       & 173           & 58            & \textbf{10}   &  & 512.0           & 143.6          & \textbf{23.4}    \\
			rajat31             & 2840          & 2767          & \textbf{1500} &  & 9590.1          & 9456.2         & \textbf{5166.6}  \\  \hline \vspace{-0.2cm}\\
			\textbf{Number of bests}& \textbf{11}	& \textbf{8} & \textbf{58} & & \textbf{9}&  \textbf{11}& \textbf{56} \\ 
			\hline
		\end{tabular}
	}
\end{table}

\begin{figure}[tbhp]
	\vspace{-0.3cm}
	\begin{center} 
		\subfloat [] {  	
			\includegraphics[scale=.47]{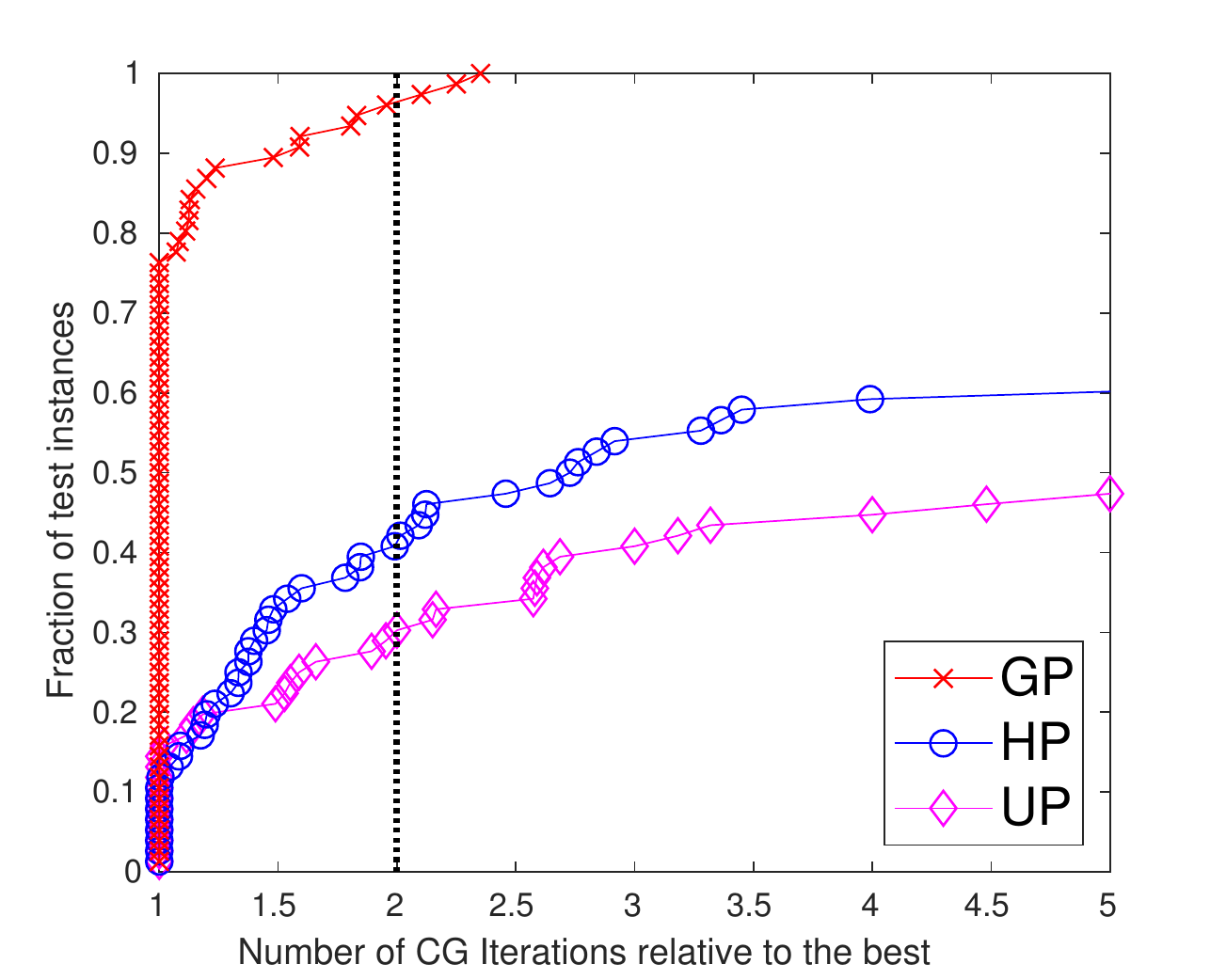}
			\label{perf1}	
		} 
		\subfloat [] {
			\includegraphics[scale=.47]{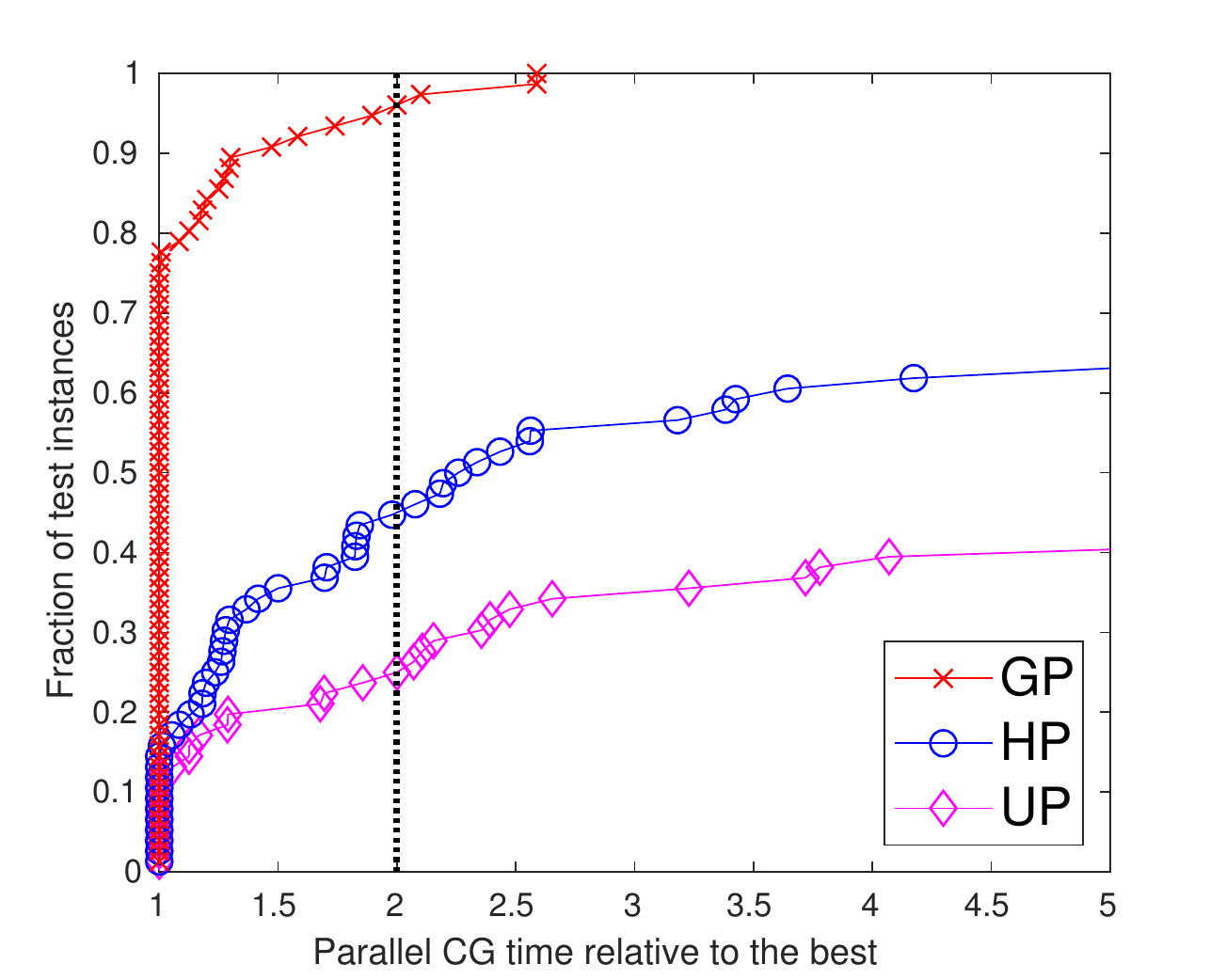}
			\label{perf2}	
		}
	\end{center}  
	\caption{ Performance profiles for (a) convergence rate and (b) parallel CG time.}
	\label{perf0}
	\vspace{-0.4cm}
\end{figure}

\Cref{perf0} shows the performance profiles \cite{dolan2002benchmarking} of $76$ matrix instances which compare multiple methods over a set of test instances.
A performance profile is used to compare multiple methods with respect to the best performing method for each test instance 
and report the fraction of the test instances in which performance is within a factor of that of the best method.
For example, in~\cref{perf1} a point (abscissa=$2$, ordinate=$0.40$) on the performance curve of a method refers to the fact that the performance of the respective method is no worse than that of the best method by a factor of 2 in approximately $40\%$ of the instances.
If a method is closer to the top-left corner, then it achieves a better performance.

In~\cref{perf1,perf2}, we show the performance profiles in terms of the number of CG iterations and parallel CG times, respectively.
As seen in~\cref{perf1}, the number of CG iterations required by GP for convergence does not exceed that of the best method by a factor of 2 in approximately 95\% of the instances, whereas HP and UP
achieve the same relative performance compared to the best method in approximately 42\% and 30\% of the instances, respectively.
As seen in~\cref{perf2}, the CG time using GP is not slower than that of the best method by a factor of 2 in approximately 95\%
of instances. 
Whereas HP and UP achieve the same relative performance compared to the best method in approximately 45\% and 23\% of the instances, respectively.
\Cref{perf1,perf2} show that the CG time is directly proportional to the number of CG iterations  as expected.

It is clear that GP mainly aims at improving the convergence rate, whereas HP mainly aims at reducing total communication volume.
We made additional measurements in order to discuss this trade off between these two methods. 
Due to the lack of space, here we only summarize the average results using 32 processors.
Although GP incurs 44\% more total communication volume than HP per iteration this results in only 6\% increase the per-iteration execution time.
Here, the per-iteration execution time for a given instance shows the overall parallelization efficiency attained by the respective partitioning method irregardless of its convergence performance.
This experimental finding can be attributed to several other factors affecting the communication overhead in addition to the total communication volume as well as the fact that per-iteration execution time is dominated by the computational cost of local solution of the linear systems via a direct solver in block Cimmino.
On average, although GP incurs only 6\% more per-iteration time than HP, GP requires 59\% smaller number of iterations for convergence than HP.
This explains the significant overall performance improvement achieved by GP against HP.

\subsection{Preprocessing Overhead and Amortization}
In this subsection, we analyze the relative preprocessing overhead of the three methods UP, HP and GP in order to find out whether the intelligent partitioning methods HP and GP are amortized.
For all of the three block-row partitioning methods, the preprocessing overhead includes  matrix scaling, block-row partitioning, creating the submatrices corresponding to these block rows and distributing these submatrices among processors.
Recall that the partitioning for UP is straightforward and hence it incurs only negligible additional cost to the preprocessing overhead.
On the other hand, intelligent partitioning algorithms utilized in HP and GP incur considerable amount of cost to the overall preprocessing overhead.
Furthermore, for GP,  construction of row inner-product graph also incurs significant amount of cost.

In~\cref{pretable}, we display total preprocessing time and total execution time for all three methods for each one of the 76 test instances.
In the table, total execution time is the sum of the total preprocessing and the total solution time.
Here the total solution time includes parallel factorization and parallel CG solution times.
Note that the factorization in block Cimmino algorithm in ABCD Solver is performed by parallel sparse direct solver MUMPS~\cite{MUMPS}.
The factorization, which needs to be performed only once, involves symbolic and numerical factorizations of the coefficient matrices of the augmented systems that arise in the block Cimmino algorithm.
Note that this factorization process is embarrassingly parallel since the factorization of the coefficient matrices of the augmented systems are done independently.
Recall that MUMPS is also used during the iterative solution stage, where at each iteration a linear system is solved using the factors of the augmented systems that were computed during the factorization stage.

\begin{table}[htbhp]
	\centering
	\caption{Preprocessing time and total execution time (including preprocessing) in seconds }
	\label{pretable}
	\centering
	\scriptsize\renewcommand{\arraystretch}{.9}
	{
		\begin{tabular}{lrrrrrrr}
			\hline \vspace{-0.15cm}\\
			\multirow{3}{*}{Matrix name}        & \multicolumn{3}{l}{Preprocessing time} &  & \multicolumn{3}{l}{Total execution time} \\  
			
			\cmidrule(l){2-4}  \cmidrule(l){6-8}  
			& UP   & HP    & GP    &  & UP   & HP   & GP      \\ \hline
			\vspace{-0.15cm}\\
			
			rajat26             & 0.12  & 0.53   & 0.52  &  & F               & 63.7           & \textbf{9.3}     \\
			ecl32               & 0.20  & 0.98   & 0.54  &  & 246.7           & 39.3           & \textbf{11.8}    \\
			2D\_54019\_highK    & 0.22  & 0.69   & 0.50  &  & 1.5             & 1.4            & \textbf{1.3}     \\
			bayer01             & 0.14  & 0.41   & 0.39  &  & 64.4            & 9.7            & \textbf{4.0}     \\
			TSOPF\_RS\_b39\_c30 & 0.47  & 0.93   & 0.87  &  & 26.1            & \textbf{6.2}   & 9.4              \\
			venkat01            & 0.75  & 2.02   & 1.44  &  & 3.0             & 3.4            & \textbf{2.8}     \\
			venkat25            & 0.67  & 2.05   & 1.55  &  & 32.0            & 16.9           & \textbf{16.7}    \\
			venkat50            & 0.76  & 2.07   & 1.46  &  & 57.3            & 26.0           & \textbf{25.8}    \\
			laminar\_duct3D     & 1.30  & 5.89   & 4.63  &  & 34.3            & 53.8           & \textbf{33.8}    \\
			lhr71c              & 0.60  & 1.82   & 1.68  &  & F               & 197.0          & \textbf{140.3}   \\
			shyy161             & 0.16  & 0.50   & 0.32  &  & \textbf{1.2}    & 1.7            & 1.5              \\
			circuit\_4          & 0.16  & 0.79   & 1.13  &  & F               & 16.6           & \textbf{13.1}    \\
			epb3                & 0.21  & 0.70   & 0.44  &  & 88.3            & 102.1          & \textbf{78.7}    \\
			poisson3Db          & 1.10  & 6.02   & 5.08  &  & 723.0           & 60.2           & \textbf{48.6}    \\
			rajat20             & 0.31  & 2.05   & 2.05  &  & 83.0            & 44.1           & \textbf{37.8}    \\
			rajat25             & 0.28  & 2.10   & 2.01  &  & 122.1           & 66.5           & \textbf{49.9}    \\
			rajat28             & 0.30  & 2.02   & 2.06  &  & 87.9            & 25.8           & \textbf{24.4}    \\
			LeGresley\_87936    & 0.35  & 0.88   & 0.64  &  & F               & 294.7          & \textbf{123.5}   \\
			rajat16             & 0.24  & 1.97   & 1.82  &  & 836.5           & 71.8           & \textbf{22.9}    \\
			ASIC\_100ks         & 0.37  & 1.50   & 1.32  &  & F               & \textbf{3.1}   & 4.1              \\
			ASIC\_100k          & 0.51  & 2.24   & 1.62  &  & 27.7            & 99.5           & \textbf{18.2}    \\
			matrix\_9           & 0.57  & 2.19   & 1.81  &  & \textbf{275.0}  & F              & 753.9            \\
			hcircuit            & 0.27  & 0.73   & 0.81  &  & 593.4           & 137.8          & \textbf{34.2}    \\
			lung2               & 0.24  & 0.52   & 0.46  &  & \textbf{1.6}    & 1.8            & 1.9              \\
			rajat23             & 0.31  & 1.12   & 1.77  &  & F               & 467.9          & \textbf{37.9}    \\
			Baumann             & 0.39  & 1.39   & 0.81  &  & \textbf{49.4}   & 89.5           & 92.5             \\
			barrier2-1          & 1.10  & 4.40   & 3.58  &  & F               & F              & \textbf{150.2}   \\
			barrier2-2          & 1.10  & 4.40   & 3.56  &  & F               & F              & \textbf{109.0}   \\
			barrier2-3          & 1.10  & 4.42   & 3.45  &  & F               & F              & \textbf{121.0}   \\
			barrier2-4          & 1.10  & 4.35   & 3.58  &  & F               & F              & \textbf{125.8}   \\
			barrier2-10         & 1.20  & 4.42   & 3.61  &  & F               & F              & \textbf{147.2}   \\
			barrier2-11         & 1.20  & 4.40   & 3.66  &  & F               & F              & \textbf{136.3}   \\
			barrier2-12         & 1.10  & 4.47   & 3.61  &  & F               & F              & \textbf{128.0}   \\
			barrier2-9          & 1.10  & 4.32   & 3.54  &  & F               & F              & \textbf{145.7}   \\
			torso2              & 0.46  & 1.17   & 0.86  &  & \textbf{1.8}    & 2.6            & 2.4              \\
			torso1              & 3.10  & 35.19  & 7.59  &  & F               & \textbf{421.5} & 983.1            \\
			dc1                 & 0.37  & 1.82   & 1.49  &  & 287.7           & 752.9          & \textbf{51.3}    \\
			dc2                 & 0.35  & 1.77   & 1.44  &  & 204.0           & 514.1          & \textbf{30.3}    \\
			dc3                 & 0.37  & 1.80   & 1.51  &  & 828.1           & 1233.9         & \textbf{48.5}    \\
			trans4              & 0.39  & 2.01   & 1.54  &  & 129.5           & 462.0          & \textbf{9.9}     \\
			trans5              & 0.40  & 2.05   & 1.42  &  & 377.3           & 1705.9         & \textbf{12.4}    \\
			matrix-new\_3       & 0.42  & 1.90   & 1.37  &  & F               & 715.0          & \textbf{524.1}   \\
			cage12              & 0.97  & 5.47   & 3.99  &  & \textbf{23.9}   & 64.2           & 33.6             \\
			FEM\_3D\_thermal2   & 1.40  & 4.90   & 3.75  &  & \textbf{8.3}    & 11.6  & 10.3             \\
			para-4              & 1.40  & 6.62   & 4.99  &  & 1605.4          & F              & \textbf{478.3}   \\
			para-10             & 1.00  & 4.45   & 3.07  &  & F               & F              & \textbf{673.6}   \\
			para-5              & 0.91  & 4.82   & 3.20  &  & F               & F              & \textbf{619.7}   \\
			para-6              & 1.10  & 4.57   & 3.13  &  & F               & F              & \textbf{610.2}   \\
			para-7              & 1.00  & 4.62   & 3.20  &  & F               & F              & \textbf{552.0}   \\
			para-8              & 1.10  & 4.57   & 3.07  &  & F               & F              & \textbf{608.8}   \\
			para-9              & 1.00  & 4.72   & 3.23  &  & F               & F              & \textbf{691.2}   \\
			crashbasis          & 0.84  & 2.50   & 2.01  &  & 13.2            & 5.2            & \textbf{5.1}     \\
			majorbasis          & 0.72  & 2.40   & 2.11  &  & 9.6             & 5.1            & \textbf{5.0}     \\
			ohne2               & 2.40  & 12.24  & 9.27  &  & F               & 2128.6         & \textbf{843.9}   \\
			hvdc2               & 0.79  & 1.82   & 1.34  &  & F               & 449.4          & \textbf{264.2}   \\
			shar\_te2-b3        & 0.16  & 2.97   & 3.10  &  & \textbf{24.2}   & 35.2           & 35.8             \\
			stomach             & 1.20  & 5.12   & 3.21  &  & \textbf{5.1}    & 9.5            & 7.8              \\
			torso3              & 1.90  & 9.27   & 4.90  &  & \textbf{14.1}   & 25.5           & 17.4             \\
			ASIC\_320ks         & 1.10  & 3.02   & 3.02  &  & F               & 8.7            & \textbf{5.2}     \\
			ASIC\_320k          & 1.30  & 12.45  & 5.34  &  & 61.9            & 64.9           & \textbf{20.4}    \\
			ML\_Laplace         & 9.30  & 58.70  & 29.05 &  & \textbf{3913.4} & 4269.0         & 4288.5           \\
			RM07R               & 14.00 & 110.00 & 69.14 &  & F               & F              & \textbf{10727.5} \\
			language            & 1.20  & 4.50   & 4.32  &  & 338.5           & 204.8          & \textbf{160.3}   \\
			CoupCons3D          & 6.40  & 32.72  & 17.66 &  & 109.5           & 103.0          & \textbf{84.5}    \\
			largebasis          & 2.20  & 6.37   & 4.61  &  & 553.1           & 182.9          & \textbf{98.8}    \\
			cage13              & 4.20  & 30.25  & 17.10 &  & \textbf{78.6}   & 301.1          & 217.3            \\
			rajat30             & 2.90  & 19.71  & 14.99 &  & 326.9           & 350.0          & \textbf{127.5}   \\
			ASIC\_680k          & 2.60  & 16.21  & 7.79  &  & \textbf{19.3}   & 228.2          & 23.8             \\
			atmosmodd           & 5.60  & 36.38  & 11.97 &  & \textbf{881.9}  & 2165.8         & 1464.6           \\
			atmosmodj           & 6.80  & 36.17  & 12.29 &  & \textbf{1006.7}          & 2372.9         & 1535.1  \\
			Hamrle3             & 4.70  & 20.00  & 9.92  &  & F               & \textbf{2406.8}         & 2469.9  \\
			atmosmodl           & 7.20  & 45.18  & 14.71 &  & 1815.8          & 1067.0         & \textbf{585.5}   \\
			atmosmodm           & 8.00  & 45.68  & 14.07 &  & 1747.3          & 934.9          & \textbf{333.2}   \\
			memchip             & 10.00 & 43.12  & 20.15 &  & 7476.3          & 1640.0         & \textbf{893.4}   \\
			circuit5M\_dc       & 16.00 & 58.06  & 25.95 &  & 540.3           & 209.7 & \textbf{64.5}    \\
			rajat31             & 25.00 & 74.36  & 37.53 &  & 9634.0          & 9539.1         & \textbf{5220.2}  \\
			\hline  \vspace{-0.2cm}\\
			&  \multicolumn{3}{c}{Geometric means} & &  \multicolumn{3}{c}{   Number of bests} \\
			& \textbf{0.9}	& \textbf{4.2} & \textbf{	2.9	} & & \textbf{15}	& \textbf{4} & \textbf{57} \\
			\hline
		\end{tabular}
	}
\end{table}

As seen in~\cref{pretable}, in terms of the preprocessing time, UP is the clear winner in all of the $76$ instances as expected, whereas GP incurs much less preprocessing time than HP in all except 7 instances.
As also seen in~\cref{pretable}, comparing all three methods, GP achieves the smallest total execution time in 57 instances, whereas UP and HP respectively achieve the smallest total
execution time in only 15 and 4 instances.

In the relative comparison of GP and UP, GP incurs only $209\%$ more preprocessing time than UP, on average, which leads GP to achieve less total execution time than UP in 61 out of 76 instances.
In other words, compared to UP, the sequential implementation of GP amortizes its preprocessing cost for those 61 instances by reducing the number of CG iterations sufficiently.
That is, GP amortizes its preprocessing cost for the solution of those 61 instances even if we solve only one linear system with a single right-hand side vector.
Note that in many applications in which sparse linear systems arise,  the solution of consecutive linear systems are required where the coefficient matrix remains the same but the right-hand side vectors change.
In such applications, the amortization performance of the proposed GP method will further improve. 
For example, the amortization performance of GP will improve from 61 to 64 instances for solving only two consecutive linear systems.

\begin{figure}[htbhp]
	\begin{center} 
		\vspace{-0.4cm}
		\subfloat [] 
		{
			\includegraphics[scale=.47]{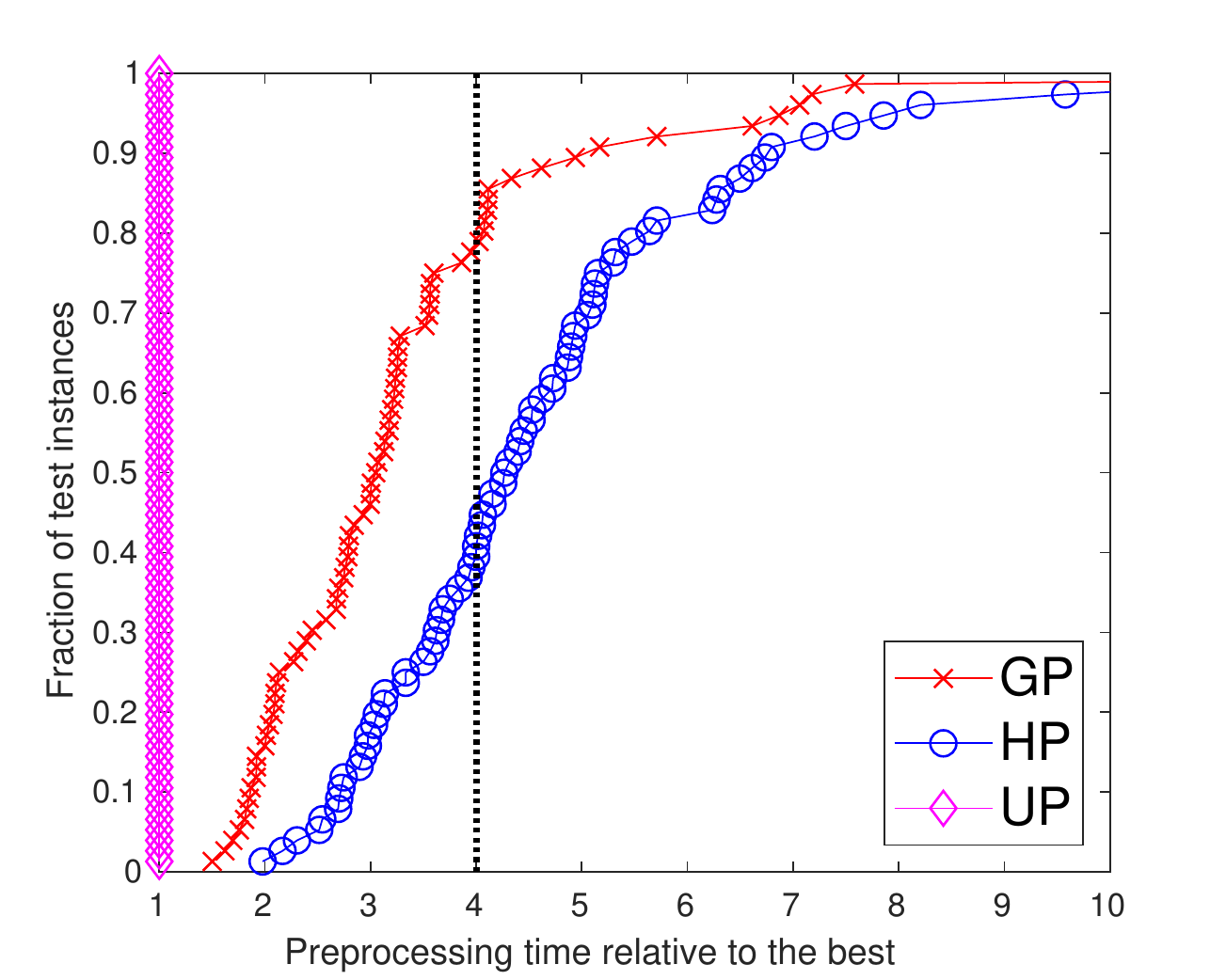}
			\label{s_part2}	
		}	
		\subfloat [] {
			\includegraphics[scale=.47]{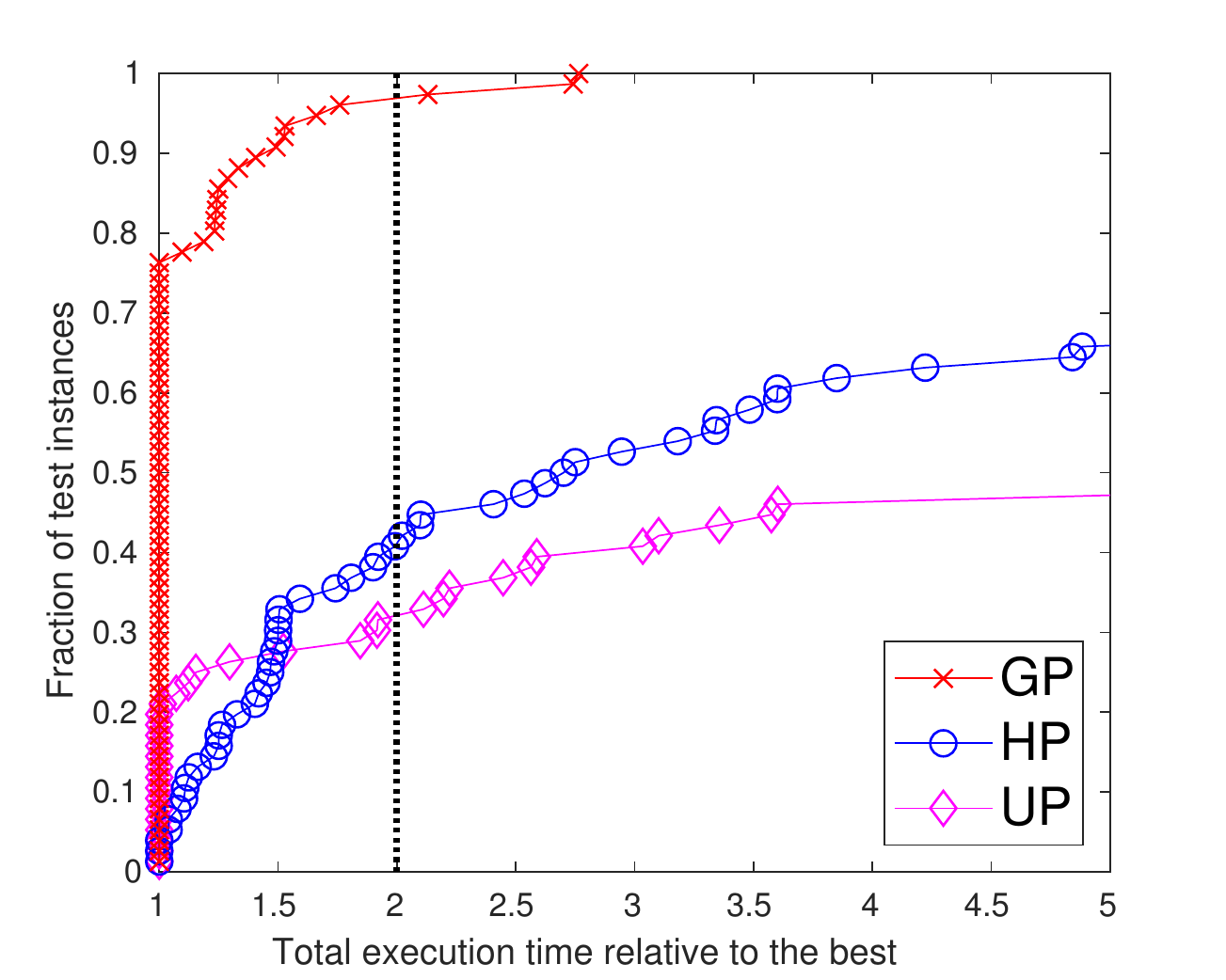}
			\label{total}	
		} 	
	\end{center}  
	\caption{Performance profiles for (a) preprocessing time and (b) total execution time for GP, HP and UP methods.}
	\label{perf}
	\vspace{-0.4cm}
\end{figure}

In~\cref{s_part2,total}, we show the performance profiles in terms of the preprocessing time and the total execution time, respectively.
As seen in~\cref{s_part2}, the preprocessing overhead incurred by GP remains no worse than that of the best method UP by a factor of 4 in approximately $77\%$ of instances.
As seen in~\cref{total}, the total execution time attained by GP does not exceed that of the best method by a factor of 2 in $96\%$ of the instances. 
On the other hand, HP and UP achieves the same relative performance compared to the best method only in approximately $44\%$ and $33\%$ of the instances, respectively.

\section{Conclusion and future work}
\label{section4}

In this paper, we proposed a novel partitioning method in order to improve the CG accelerated block Cimmino algorithm.
The proposed partitioning method takes the numerical orthogonality between block rows of the coefficient matrix into account.   
The experiments on a large set of real world systems show that the  proposed method improves the convergence rate of the CG accelerated block Cimmino compared to the state-of-the-art hypergraph partitioning method. 
Moreover, it requires not only less preprocessing time and fewer number of CG iterations, but also much less total execution time than the hypergraph partitioning method.   

As a future work, we consider two issues:
further reducing the number of iterations through preconditioning and reducing the preprocessing overhead through parallelization.

Even though the $H$ matrix is not available explicitly, it could be still possible to obtain a preconditioner, further reducing the required number of CG iterations.  One viable option could be using sparse approximate inverse type preconditioner where the coefficient matrix does not need to be available explicitly. This approach could be viable especially when consecutive linear systems are needed to be solved with the same coefficient matrix.        

The proposed method involves two computational stages, namely constructing row inner-product graph via computing SpGEMM operation and partitioning this graph.
For parallelizing the first stage, parallel SpGEMM~\cite{akbudak2014simultaneous,kadirSpGEMM,azad2016exploiting} operation could be used to construct local subgraphs on each processor.
For parallelizing the second stage, a parallel graph partitioning tool ParMETIS~\cite{karypis1998parallel} could be used.
In each processor, the local subgraphs generated in parallel in the first stage could be used as input for ParMETIS.

\appendix

\section{Graph and Graph Partitioning}
\label{appendix}
An undirected graph $\mathcal{G} = (\mathcal{V}, \mathcal{E})$ consists of a set of vertices $\mathcal{V}$ and a set of edges $\mathcal{E}$.
Each edge $(v_i,v_j) \in \mathcal{E}$ connects a pair of distinct vertices $v_i$ and $v_j$.
Each vertex $v_i \in \mathcal{V}$ can be assigned a weight shown as $w(v_i)$ and each edge $(v_i,v_j)  \in \mathcal{E}$ can be assigned a cost shown as $cost(v_i,v_j)$.

$\Pi = \{ \mathcal{V}_1,\mathcal{V}_2,\ldots,\mathcal{V}_K \}$ is defined as a K-way vertex partition of $\mathcal{G}$  if parts are mutually disjoint and exhaustive.
An edge $(v_i,v_j)$ is said to be cut if the vertices $v_i$ and $v_j$ belong  to different vertex parts  and uncut otherwise.
The set of cut edges of a partition $\Pi$ is denoted as $\mathcal{E}_{\rm{cut}}$.
In a given partition $\Pi$ of $\mathcal{G}$, the weight $W_k$ of a part $\mathcal{V}_k$ is defined as the sum of the weights of the vertices in $\mathcal{V}_k$, i.e.,
\begin{equation}
\label{Wk}
W_k \triangleq \sum\limits_{v_i \in \mathcal{V}_k} w(v_i).
\end{equation}

\noindent In the graph partitioning problem, the partitioning constraint is to maintain a given  balance condition on the part weights, i.e.,
\begin{equation}
\label{partbal}
W_{k} \leq W_{avg}(1+\epsilon) \text{ for } k =1,2\ldots,K,\text{ where
} W_{avg} = \sum\limits_{v_i \in \mathcal{V}}w(v_i) /K. 
\end{equation}
Here $\epsilon$ is the predefined maximum imbalance ratio. 
The partitioning objective is to minimize the cutsize defined as the sum of the costs of the cut edges, i.e.,
\begin{equation}
\label{partcut}
{\rm{cutsize}}(\Pi) \triangleq \sum\limits_{ (v_i,v_j) \in \mathcal{E}_{cut}} cost(v_i,v_j).
\end{equation}

\section*{Acknowledgments}
This work was done when the first author was at Bilkent University, Ankara, Turkey and was revised when the first author was at IRIT-CNRS, Toulouse, France.

\bibliographystyle{siamplain} 
\bibliography{cites}
\end{document}